\documentclass{amsart}
\usepackage{lamsarrow}
\usepackage{pb-diagram}
\usepackage{amssymb}
\usepackage{pb-lams}
\numberwithin{equation}{subsection}
\newtheorem{theorem}[equation]{Theorem}
\newtheorem{lemma}[equation]{Lemma}
\newtheorem{corollary}[equation]{Corollary}
\newtheorem{proposition}[equation]{Proposition}
\newtheorem{definition}[equation]{Definition}
\newtheorem{remark}[equation]{Remark}
\newtheorem{example}[equation]{Example}

\DeclareMathOperator{\Hom}{Hom}
\DeclareMathOperator{\RHom}{RHom}
\DeclareMathOperator{\End}{End}

\DeclareMathOperator{\map}{map}

\DeclareMathOperator{\sh}{sh}
\DeclareMathOperator{\Ho}{Ho}
\DeclareMathOperator{\HoC}{Ho(\C)}
\DeclareMathOperator{\hocolim}{hocolim}

\DeclareMathOperator{\Ev}{Ev}

\DeclareMathOperator{\f}{f}

\newcommand{\mN}{{\mathbb N}}
\newcommand{\mS}{{\mathbb S}}

\newcommand{\mZ}{{\mathbb Z}}
\newcommand{\uA}{\underline{A}}
\newcommand{\A}{{\mathcal A}}
\newcommand{\C}{{\mathcal C}}
\newcommand{\D}{{\mathcal D}}
\newcommand{\E}{\mathcal E}
\newcommand{\G}{\mathcal G}
\newcommand{\R}{\mathcal R}
\newcommand{\T}{\mathcal T}
\newcommand{\U}{\mathcal U}

\newcommand{\Hc}{\mathcal H}
\newcommand{\Oc}{\mathcal O}
\newcommand{\Pc}{\mathcal P}
\newcommand{\EP}{\mathcal E(\mathcal P)}
\newcommand{\EG}{\mathcal E(\mathcal G)}

\newcommand{\Spc}{{\mathcal S}p}
\newcommand{\spec}{Sp^{\Sigma}}
\newcommand{\specC}{Sp(\C)}

\newcommand{\sset}{\mathcal{S}_*}
\newcommand{\Omod}{$\Oc$-module}
\newcommand{\Omods}{$\Oc$-modules}
\newcommand{\boxprod}{\mathbin{\square }}

\newcommand{\iso}{\cong}
\newcommand{\rscript}[1]{^{\mbox{\scriptsize {#1}}}}
\newcommand{\rtiny}[1]{^{\mbox{\tiny {#1}}}}
\newcommand{\Sigf}{\Sigma^{\infty}_{\f}}
\newcommand{\sm}{\wedge}

\newcommand{\tensor}{\otimes}

\newcommand{\widebar}{\overline}

\renewcommand{\to}{\longrightarrow}
\newcommand{\varrow}[1]{\hbox to #1{\rightarrowfill}}
\newcommand{\varl}[2]{\stackrel{#2}{\hbox to #1{\leftarrowfill}}}
\newcommand{\varr}[2]{\stackrel{#2}{\hbox to #1{\rightarrowfill}}}
\newcommand{\parallelarrows}[1]{\begin{array}{c} {\hbox to
#1{\rightarrowfill}}  \vspace{-0.35cm} \\ {\hbox to
#1{\rightarrowfill}} \end{array}}

\begin{document}

\title[Stable model categories]{Classification of stable model  
categories 
}
\date{\today; 1991 AMS Math.\ Subj.\ Class.: 55U35, 55P42}
\author{Stefan Schwede}
\thanks{Research supported by a BASF-Forschungsstipendium der
Studienstiftung des deutschen Volkes}
\address{Fakult\"at f\"ur Mathematik \\ Universit\"at Bielefeld\\
33615 Bielefeld, Germany}
\email{schwede@mathematik.uni-bielefeld.de}
\author{Brooke Shipley}
\thanks{Research partially supported by an NSF Postdoctoral Fellowship}
\address{Department of Mathematics \\ Purdue University\\  
W. Lafayette, IN 47907 
\\ USA}
\email{bshipley@math.purdue.edu}
\maketitle

{\small Abstract:
A stable model category is a setting for homotopy theory
where the suspension functor is invertible.
The prototypical examples are the category of spectra in the sense of stable
homotopy theory and the category of unbounded chain complexes of modules over
a ring.
In this paper we develop methods for deciding when two stable model
categories represent `the same homotopy theory'.
We show that stable model categories with a single compact generator are
equivalent to modules over a ring spectrum.
More generally stable model categories with a set of generators are
characterized as modules over a `ring spectrum with several objects',
i.e., as spectrum valued diagram categories.
We also prove a Morita theorem which shows how equivalences between module
categories over ring spectra can be realized by smashing with a pair of
bimodules. Finally, we characterize stable model categories which
represent the derived category of a ring.  This is a slight
generalization of Rickard's work on derived equivalent rings.
We also include a proof of the model category equivalence
of modules over the Eilenberg-Mac Lane spectrum $HR$ and (unbounded) chain 
complexes of $R$-modules for a ring $R$.}

\section{Introduction}

The recent discovery of highly structured categories of spectra has opened
the way for a new wholesale use of algebra in stable homotopy theory.
In this paper we use this new algebra of spectra to characterize stable
model categories, the settings for doing stable homotopy theory, as
categories of highly structured modules.  This characterization also leads
to a Morita theory for equivalences between categories of highly
structured modules.

The motivation and techniques for this paper come from two directions,
namely stable homotopy theory and homological algebra.
Specifically, stable homotopy theory studies the classical stable homotopy
category which is the category of spectra up to homotopy.
For our purposes though, the homotopy category is inadequate because too much 
information is lost, for example the homotopy type of mapping spaces.  
Instead, we study the model category of spectra which
captures the whole stable homotopy theory.
More generally we study stable model categories, those model categories
which share the main formal property of spectra, namely that the suspension
functor is invertible up to homotopy.
We list examples of stable model categories in Section \ref{definitions}.

The algebraic part of the motivation arises as follows. A classical theorem,
due to Gabriel \cite{gabriel},
characterizes categories of modules as the cocomplete
abelian categories with a single small
projective generator; the classical Morita theory for equivalences 
between module categories (see for example \cite[\S 21, 22]{anderson-fuller})
follows from this.
Later Rickard \cite{rickard1,rickard2} developed a Morita theory 
for derived categories based on the notion of a tilting complex.
In this paper we carry this line of thought one step further.
Spectra are the homotopy theoretical generalization of abelian
groups and stable model categories are the homotopy
theoretic analogue of abelian categories
(or rather their categories of chain complexes).  
Our generalization of Gabriel's theorem develops a Morita theory 
for stable model categories.
Here the appropriate notion of a model category equivalence 
is a Quillen equivalence since these equivalences preserve the homotopy theory,
not just the homotopy category, see~\ref{sec-Quillen-equiv}.
 
We have organized our results into three groups:
\smallskip

{\bf Characterization of module categories.}
The model category of modules over a
ring spectrum has a single compact generator,
namely the free module of rank one.
But module  categories are actually characterized by this property.
To every object in  a stable model category we associate an
endomorphism ring spectrum, see Definition \ref{def-endo ring spectrum}.
We show that if there is a single
compact generator, then the given stable model category has the same homotopy
theory as the modules over the endomorphism ring spectrum of the generator
(Theorem \ref{thm-main-one}).
More generally in Theorem~\ref{thm-main-gen}, 
stable model categories with a set of compact generators
are characterized as modules over a `ring spectrum with many objects,' or
{\em spectral category}, see Definition~\ref{def-spectral cat}.
This is analogous to Freyd's generalization of 
Gabriel's theorem~\cite[5.3H p.120]{freyd}.
Examples of these characterizations are given 
in~\ref{single generator examples} and~\ref{multiple generator examples}.

\medskip

{\bf Morita theory for ring spectra.}
In the classical algebraic context Morita theory describes equivalences
between module categories in terms of bimodules, 
see e.g.\ \cite[Thm.\ 22.1]{anderson-fuller}.
In Theorem~\ref{thm-Morita} we present an analogous result 
which explains how a chain of Quillen equivalences between module categories 
over ring spectra can be replaced by a single Quillen equivalence 
given by smashing with a pair of bimodules.

\medskip
 
{\bf Generalized tilting theory.}
In \cite{rickard1, rickard2}, Rickard answered the question
of when two rings are derived equivalent, i.e., when various derived module
categories are equivalent as triangulated categories.
Basically, a derived equivalence exists if and only if
a so-called tilting complex exists. From our point of view, 
a tilting complex is a particular compact generator
for the derived category of a ring.
In Theorem \ref{thm-tilting} we obtain a generalized tilting theorem
which characterizes stable model
categories which are Quillen (or derived) equivalent to the derived category
of a ring.

\medskip

Another result which is very closely related to this characterization of
stable model categories can be found in~\cite{ss-uniqueness} where we
give necessary and sufficient conditions for when a stable model
category is Quillen equivalent to spectra, 
see also Example~\ref{single generator examples} (i).
These uniqueness results are then developed further in~\cite{ship-mon, sch-2}.
Moreover, the results in this paper form a basis for developing
an algebraic model for any rational stable model category.  This
is carried out in~\cite{ss-equiv} and applied 
in~\cite{shipley-S1, greenlees-shipley}.

In order to carry out our program it is essential to have available a
highly structured model for the category of spectra which admits a symmetric
monoidal and homotopically well behaved smash product before passing to the
homotopy category. The first examples of such categories were the
$S$-modules of \cite{ekmm} and the symmetric spectra of Jeff Smith \cite{hss};
by now several more such categories have been constructed
\cite{lydakis-simplicial,mmss}.
We work with symmetric spectra because we can replace stable model
categories by Quillen equivalent ones which are
enriched over symmetric spectra (Section \ref{symmetric objects}).
Also, symmetric spectra are reasonably easy
to define and understand and several other model categories in the
literature are already enriched over symmetric spectra.
The full strength of our viewpoint comes from combining enriched
(over symmetric spectra) category theory with the language 
of closed model categories.  
We give specific references throughout;
for general background on model categories 
see Quillen's original article~\cite{Q},
a modern introduction~\cite{DS}, or~\cite{hovey-book} for
a more complete overview.
 
We want to point out the conceptual similarities between the present paper
and the work of Keller \cite{keller-derivingDG}.
Keller uses differential graded categories to give an elegant reformulation
(and generalization) of Rickard's results on derived equivalences for rings.
Our approach is similar to Keller's, but where he considers categories 
whose hom-objects are chain complexes of abelian
groups, our categories have hom-objects which are spectra.
Keller does not use  the language of
model categories, but the `P-resolutions' of \cite[3.1]{keller-derivingDG}
are basically cofibrant-fibrant replacements.

{\bf Notation and conventions:} 
We use the symbol $\sset$ to denote the category of {\em pointed} simplicial 
sets, and we use $\spec$ for the category of symmetric spectra \cite{hss}. 
The letters $\C$ and $\D$ usually denote model categories,
most of the time assumed to be simplicial and stable.
For cofibrant or fibrant approximations of objects in a model
category we use superscripts $(-)\rscript{c}$ and $(-)\rscript{f}$.
For an object $X$ in a pointed simplicial model category we use
the notation $\Sigma X$ and $\Omega X$ for the simplicial suspension
and loop functors (i.e., the pointed tensor and cotensor
of an object $X$ with the pointed simplicial circle
$S^1=\Delta[1]/\partial\Delta[1]$); one should keep in mind
that these objects may have the `wrong' homotopy type if $X$ is
not cofibrant or fibrant respectively.
Our notation for various kinds of morphism objects is as follows:
the set of morphisms in a category $\C$ is denoted `$\hom_{\C}$';
the simplicial set of morphisms in a simplicial category
is denoted `$\map$'; we use `$\Hom$' for
the symmetric function spectrum in a spectral model category
(Definition \ref{def-spectral model cat});
square brackets `$[X,Y]^{\HoC}$' denote the abelian group of morphisms
in the homotopy category of a stable model category $\C$;
and for  objects $X$ and $Y$ in any triangulated category $\T$
we use the notation `$[X,Y]^{\T}_*$' to denote the {\em graded}
abelian group of morphisms, i.e., $[X,Y]_n^{\T} = [X[n],Y]^{\T}$
for $n\in{\mathbb Z}$ and where $X[n]$ is the $n$-fold shift of $X$.
 
We want to write the evaluation of a morphism $f$ on an element $x$ as $f(x)$. 
This determines the following conventions about actions of rings and ring
spectra: the endomorphism monoid, ring or ring spectrum $\End(X)$ acts 
on the object $X$ from the {\em left}, and it acts on the set 
(group, spectrum) $\Hom(X,Y)$  from the {\em right}. 
A module will always be a right module; 
this way the left multiplication map establishes an isomorphism between a
ring and the endomorphism ring of the free module of rank one. 
A $T$-$R$-bimodule is a $(T\rtiny{op}\tensor R)$-module 
(or a  $(T\rtiny{op}\sm R)$-module in the context of ring spectra).
  
{\bf Organization:}
In Section \ref{definitions} we recall
stable model categories and some of their properties,
as well as the notions of compactness and generators, and
we give an extensive list of examples.
In Section \ref{sec-classification theorem} we prove the 
classification theorems (Theorems \ref{thm-main-one} and \ref{thm-main-gen}).
In Section \ref{symmetric objects} we introduce the category
$\specC$ of symmetric spectra over a simplicial model category $\C$.
Under certain technical assumptions we  show in Theorem \ref{thm-specC}
that it is a stable model category with composable and homotopically
well-behaved function symmetric spectra which is Quillen equivalent to
the original stable model category $\C$.
In Definition \ref{def-endo ring spectrum} we associate to an object $P$
of a simplicial stable model category a symmetric
{\em endomorphism ring spectrum} $\End(P)$.
In Theorem \ref{main theorem spectral version} we then prove
Theorem \ref{thm-main-gen} for {\em spectral model categories} 
(Definition \ref{def-spectral model cat}), such as for example $\specC$.
This will complete the classification results.
In Section \ref{sec-Morita} we prove the Morita context
(Theorem \ref{thm-Morita}) and in Section \ref{sec-tilting} 
we prove the tilting theorem (Theorem \ref{thm-tilting}).
In two appendices we consider modules over spectral categories, the homotopy
invariance of endomorphism ring spectra and the characterization of
Eilenberg-Mac Lane spectral categories.

{\bf Acknowledgments:}
We would like to thank Greg Arone, Dan Dugger, Bill Dwyer, Mark Hovey,
Jeff Smith and Charles Rezk for inspiration
and for many helpful discussions on the subjects of this paper.

\section{Stable model categories}\label{definitions}

In this section we recall stable model categories and 
some of their properties, as well as  the notions
of compactness and generators, and we give a list of examples.

\subsection{Structure on the homotopy category} 
\label{sub-structure on homotopy category}

Recall from ~\cite[I.2]{Q} of \cite[6.1]{hovey-book}
that the homotopy category of a pointed
model category supports a suspension functor $\Sigma$ with a
right adjoint loop functor $\Omega$.

\begin{definition} \label{def-stable model category}
{\em  A {\em stable model category} is a
pointed closed model category for which
the functors $\Omega$ and $\Sigma$ on the homotopy category are
inverse equivalences.}
\end{definition}

The homotopy category of a stable model category has  
a large amount of extra structure, some of which plays a role in this paper. 
First of all, it is naturally a triangulated category 
(cf.\ \cite{verdier} or \cite[A.1]{hps}). 
A complete reference for this fact can be found in \cite[7.1.6]{hovey-book}; 
we sketch the constructions: by definition of `stable' the suspension 
functor is a self-equivalence of the
homotopy category and it defines the shift functor. 
Since every object is a two-fold suspension,
hence an abelian co-group object, the homotopy category of 
a stable model category is additive.
Furthermore, by \cite[7.1.11]{hovey-book} the cofiber sequences and 
fiber sequences of \cite[1.3]{Q} coincide
up to sign in the stable case, and they define the distinguished triangles. 
Since we required a stable model category
to have all limits and colimits, its homotopy category has infinite sums 
and products. So such a homotopy category behaves like the 
unbounded derived category of an abelian category.
This motivates thinking of a stable model category as a homotopy 
theoretic analog of an abelian category.
  
We recall the notions of compactness and generators in the context of  
triangulated categories:
 
\begin{definition}\label{defn-generator & compact}
{\em Let $\T$ be a triangulated category with infinite  
coproducts. A full triangulated subcategory of $\T$ (with shift  
and triangles induced from $\T$) is called {\em localizing} if it  
is closed under coproducts in $\T$. A set $\Pc$ of
objects of $\T$ is called a set of {\em generators} if the only  
localizing subcategory which contains the objects of $\Pc$ is 
$\T$ itself.
An object $X$ of $\T$ is {\em compact} 
(also called {\em small} of {\em finite}) 
if for any family of objects $\{A_i\}_{i\in I}$ the canonical map
\[ \bigoplus_{i\in I} \, [X, A_i]^{\T} \ \varrow{1cm} \ [X,  
\coprod_{i\in I} A_i]^{\T} \]
is an isomorphism.
Objects of a stable model category are called `generators' or `compact' 
if they are so when  considered as objects of the triangulated 
homotopy category.
}
\end{definition}

A triangulated category with infinite coproducts and a set 
of compact generators is often called {\em compactly generated}.
We avoid this terminology because of the danger of confusing it with
the terms `cofibrantly generated' and `compactly generated' 
in the context of model categories. 

\subsection{Remarks} \label{compact generator remarks} 
\renewcommand{\labelenumi}{(\roman{enumi})}
\begin{enumerate}

\item There is a convenient criterion for when a set of {\em compact} objects 
generates a triangulated category. This characterization is well known, 
but we have been unable to find a reference which proves it 
in the form we need.

\begin{lemma} \label{localizing} Let $\T$ be a triangulated  
category with infinite coproducts and let $\Pc$ be a set of compact  
objects. Then the following are equivalent:\\
(i) The set $\Pc$ generates $\T$ in the sense of Definition  
\ref{defn-generator & compact}. \\
(ii) An object $X$ of $\T$ is trivial if and only if there are  
no graded maps from objects of $\Pc$ to $X$, i.e.
$[P, X]_*=0$ for all $P \in \Pc$.
\end{lemma}
\begin{proof}
Suppose the set $\Pc$ generates $\T$ and let $X$ be  an object  
with the property that $[P,X]_*=0$ for all $P\in \Pc$. The full  
subcategory of $\T$ of objects $Y$ satisfying $[Y,X]_*=0$ is  
localizing. Since it contains the set $\Pc$, it contains all of  
$\T$. Taking $Y=X$ we see that the identity map of $X$ is  
trivial, so $X$ is trivial.

The other implication uses the existence of Bousfield localization  
functors, which in this case is a {\em finite localization}  
first considered by Miller in the context of the stable homotopy
category \cite{miller-finite}. For every set  
$\Pc$ of compact objects in a triangulated category with infinite  
coproducts there exist functors $L_{\Pc}$ (localization) and  
$C_{\Pc}$ (colocalization) and a natural distinguished triangle
\[ C_{\Pc} X \ \varrow{1cm} \ X \ \varrow{1cm} \ L_{\Pc} X  
\ \varrow{1cm} \ C_{\Pc}X[1] \]
such that $C_{\Pc} X$ lies in the localizing subcategory generated  
by $\Pc$, and such that $[P,L_{\Pc}X]_{\ast}=0$ for all $P\in \Pc$  
and $X\in\T$; one reference for this construction is in the proof
of \cite[Prop.\ 2.3.17]{hps} see also~\ref{rem-fin-loc}. 
So if we assume condition (ii)  
then for all $X$ the localization $L_{\Pc} X$ is trivial, hence $X$  
is isomorphic to the colocalization $C_{\Pc} X$ and thus contained  
in the localizing subcategory generated by $\Pc$.
\end{proof}

\item Our terminology for `generators' is different from the use of the
term in  category theory; 
generators in our sense are sometimes called {\em weak generators} elsewhere. 
By Lemma \ref{localizing}, a set of generators detects if {\em objects} 
are trivial (or equivalently if maps in $\T$ are isomorphisms). 
This notion has to  be distinguished from that of a categorical generator 
which detects if {\em maps} are trivial. 
For example, the sphere spectra are a set of generators 
(in the sense of Definition \ref{defn-generator & compact}) for
 the stable homotopy category of spectra. Freyd's generating hypothesis 
conjectures that the spheres are a set of categorical generators for 
the stable homotopy category of {\em finite}  spectra. 
It is unknown to this day whether the generating hypothesis is true or false.

\item An object of a triangulated category is compact if and only if  
its shifts (suspension and loop objects) are. Any finite coproduct  
or direct summand of compact objects is again compact.
Compact objects are closed under extensions: if two objects in a
distinguished triangle are compact, then so is the third one. In other  
words, the full subcategory of compact objects in a triangulated  
category is {\em thick}.
There are non-trivial triangulated categories in which only the  
zero object is compact. Examples with underlying stable model  
categories arise for example as suitable Bousfield localizations of  
the category of spectra, see \cite[Cor.\ B.13]{hovey-strickland-kn}.

\item If a triangulated category has a set of generators, then the  
coproduct of all of them is a single generator.
However, infinite coproducts of compact objects are in general not compact.
So the property of having a single compact generator is something special.
In fact we see in Theorem \ref{thm-main-one} below that this  
condition characterizes the module
categories over ring spectra among the stable model categories.
If generators exist, they are far from being unique.

\item In the following we often consider stable
model categories which are cofibrantly generated. Hovey
 has shown \cite[Thm.\ 7.3.1]{hovey-book} that a cofibrantly generated model
category always has a set of generators in the sense of Definition  
\ref{defn-generator & compact} (the cofibers of any set of generating
cofibrations will do).
So having generators is not an extra condition in the situation we consider,
although these generators may not be compact.  
See~\cite[Cor.\ 7.4.4]{hovey-book} for
conditions that guarantee a set of compact generators.
\end{enumerate}

\subsection{Examples}
\label{stable examples} 
\renewcommand{\labelenumi}{(\roman{enumi})}

\begin{enumerate}
\item {\bf Spectra.} As we mentioned in the introduction, one of  
our main motivating examples is the category of spectra in the sense  
of stable homotopy theory. The sphere spectrum is a compact generator.
Many model categories of spectra have  
been constructed, for example by
Bousfield and Friedlander~\cite{BF}; 
Robinson \cite[`spectral sheaves']{robinson-spectra}; 
Jardine \cite[`$n$-fold spectra']{jardine-etale};  
Elmendorf, Kriz, Mandell and May \cite[`coordinate free spectra',
`$\mathbb L$-spectra', `$S$-modules']{ekmm}; 
Hovey, Shipley and Smith \cite[`symmetric spectra']{hss}; 
Lydakis \cite[`simplicial functors']{lydakis-simplicial}; 
Mandell, May, Schwede and Shipley \cite[`orthogonal spectra', 
`$\mathcal W$-spaces']{mmss}.

\item {\bf Modules over ring spectra.} Modules over an  
$S$-algebra~\cite[VII.1]{ekmm} or modules over a symmetric ring spectrum  
\cite[5.4.2]{hss} form proper, cofibrantly generated, simplicial, 
stable model categories, see also~\cite[Sec.\ 12]{mmss}.
In each case a module is compact if and only if it is weakly equivalent 
to a retract of a finite cell module. 
The free module of rank one is a compact generator. 
More generally there are stable model categories of modules over 
`symmetric ring spectra with several objects', or {\em spectral categories}, 
see Definition \ref{def-spectral cat} and Theorem \ref{O-modules}.

\item {\bf Equivariant stable homotopy theory.} If $G$ is a compact  
Lie group, there is a category of $G$-equivariant coordinate free  
spectra~\cite{lms} which is a stable model category. 
Modern versions of this model category are the $G$-equivariant orthogonal
spectra of~\cite{mm} and $G$-equivariant $S$-modules of~\cite{ekmm}.
In this case  the equivariant suspension spectra of the coset spaces $G/H_+$ 
for all closed subgroups $H\subseteq G$ form a set of compact generators.
This equivariant model category is taken up again in 
Examples \ref{multiple generator examples} (i) and 
\ref{ex-rat equivariant via tilting}.

\item {\bf Presheaves of spectra.}  For every Grothendieck site Jardine 
\cite{jardine-stable_presheaves}
constructs a proper, simplicial, stable model category of presheaves of 
Bousfield-Friedlander type spectra; the  weak equivalences are the maps 
which induce isomorphisms of the associated sheaves of stable homotopy groups.
For a general site these stable model categories do not seem to have a 
set of compact generators.

\item {\bf The stabilization of a model category.}
In principle every pointed model category
should give rise to an associated stable model category by  
`inverting' the suspension functor, i.e., by passage to internal  
spectra.
This has been carried out for certain simplicial model categories  
in \cite{sch-cotangent} and \cite{hovey-stabilization}. The  
construction of symmetric spectra over a model category (see Section  
\ref{symmetric objects}) is another approach to stabilization.

\item {\bf Bousfield localization.} 
Following Bousfield~\cite{bousfield-localization-spaces}, localized model 
structures for modules over an $S$-algebra are constructed in 
\cite[VIII 1.1]{ekmm}. Hirschhorn \cite{hirschhorn-book}
shows that under quite general hypotheses the localization 
of a model category is again a model category. 
The localization of a stable model category is stable and localization
preserves generators. 
Compactness need not be preserved, see Example  
\ref{single generator examples} (iii).

\item {\bf Motivic stable homotopy.} 
In \cite{morel-voevodsky, voevodsky-icm} Morel and Voevodsky introduced the 
${\mathbb A}^1$-local model category structure for schemes over a base. 
An associated stable homotopy category of ${\mathbb A}^1$-local 
$T$-spectra (where $T = {\mathbb A}^1/({\mathbb A}^1-0)$ is the `Tate-sphere') 
is an important tool in Voevodsky's proof of the Milnor conjecture 
\cite{voevodsky-milcon}. This stable homotopy category arises
from a stable model category with a set of compact generators, 
see Example \ref{multiple generator examples} (ii) for more details.
\end{enumerate}

\subsection{Examples: abelian stable model categories} \label{abelian examples} 
\renewcommand{\labelenumi}{(\roman{enumi})}

Some examples of stable model categories are
`algebraic', i.e., the model category is also an abelian category. 
Most of the time the objects consist of chain complexes
in some abelian category and depending on the choice of weak equivalences
one gets a kind of derived category as the homotopy category.
A different kind of example is formed by the stable module categories
of Frobenius rings. 

For algebraic examples as the ones below, our results are
essentially covered by Keller's paper \cite{keller-derivingDG},
although Keller does not use the language of model categories.
Also there is no need to consider spectra when dealing with 
abelian model categories: 
the second author shows \cite{shipley-abelian} that every
cofibrantly generated, proper, abelian stable model category is Quillen 
equivalent to a DG-model category, i.e., a model category enriched over 
chain complexes of abelian groups.

\begin{enumerate}
\item \label{ex-modules} {\bf Complexes of modules.} The category of  
unbounded chain complexes of left modules over a ring supports 
a model category structure with weak equivalences the quasi-isomorphisms 
and with fibrations the epimorphisms \cite[Thm.\ 2.3.11]{hovey-book}
(this is called the {\em projective} model structure).
Hence the associated homotopy category is the unbounded derived category 
of the ring. 
A chain complex of modules over a ring is compact if and only if it is  
quasi-isomorphic to a bounded complex of finitely generated  
projective modules \cite[Prop.\ 6.4]{boek-nee}.
We show in Theorem \ref{thm-chains and EM} that the model category 
of unbounded chain complexes of $A$-modules is Quillen equivalent 
to the category of modules over the symmetric  
Eilenberg-Mac Lane ring spectrum for $A$. 
This example can be generalized in at least two directions: 
one can consider model categories of chain complexes in an abelian category
with enough projectives (see e.g.\ \cite[2.2]{Christensen-Hovey-relative}
for a very general construction under mild smallness assumptions).
On the other hand one can consider model categories of differential 
graded modules over a differential graded algebra, or even a 
`DGA with many objects',  alias DG-categories~\cite{keller-derivingDG}.

\item {\bf Relative homological algebra.} 
In \cite{Christensen-Hovey-relative},
Christensen and Hovey introduces model category structures for chain complexes
over an abelian category based on a {\em projective class}.
In the special case where the abelian category is modules over some ring
and the projective class consists of all summands of free modules this
recovers the (projective) model category structure of the previous example.
Another special case of interest is the {\em pure derived category} of a ring.
Here the projective class consists of all summands of (possibly infinite) 
sums of finitely generated modules, 
see also Example \ref{ex-pure via tilting}. 

\item {\bf Homotopy categories of abelian categories.} 
For any abelian category $\A$, there is a stable model category structure
on the category of unbounded chain complexes in $\A$ with 
the {\em chain homotopy equivalences} as weak equivalences,
see e.g., \cite[Ex.~3.4]{Christensen-Hovey-relative}. 
The associated homotopy category is usually denoted $K(\A)$.
Such triangulated homotopy categories tend not to have a set of small
generators; for example, Neeman~\cite[E.3.2]{neeman-triangle book} 
shows that the homotopy category of chain complexes of abelian groups  
$K(\mZ)$ does not have a set of generators whatsoever.

\item {\bf Quasi-coherent sheaves.} For a nice enough scheme $X$
the derived category of quasi-coherent sheaves $\D(qc/X)$ arises 
from a stable model category and has a set of compact generators.
More precisely, if $X$ is quasi-compact and quasi-separated, then the so-called
{\em injective} model structure exists.
The objects of the model category are unbounded complexes of quasi-coherent
sheaves of $\Oc_X$-modules, 
the weak equivalences are the quasi-isomorphisms and the cofibrations
are the injections \cite[Cor.\ 2.3 (b)]{hovey-sheaves}.
If $X$ is separated, then the compact objects of the derived category 
are precisely the {\em perfect complexes}, 
i.e., the complexes which are locally quasi-isomorphic 
to a bounded complex of vector bundles
\cite[2.3, 2.5]{neeman-grothendieck}.
If $X$ also admits an ample family of line bundles 
$\{{\mathcal L}_{\alpha}\}_{\alpha\in A}$, then the set of line bundles
$\{{\mathcal L}_{\alpha}^{\tensor m} \, | \, \alpha\in A, m\in \mZ\}$, 
considered as complexes concentrated in dimension zero, 
generates the derived category $\D(qc/X)$, 
see \cite[1.11]{neeman-grothendieck}.
This class of examples contains the derived category of a ring 
as a special case, but the injective model structure is different from the one
mentioned in (i).
Hovey \cite[Thm.\ 2.2]{hovey-sheaves} generalizes the injective 
model structure to abelian Grothendieck categories.

\item {\bf The stable module category of a  Frobenius ring.} 
A Frobenius ring is defined by the property that  
the classes of projective and injective modules coincide. Important  
examples are finite dimensional Hopf-algebras over a field and in  
particular group algebras of finite groups.
The {\em stable module category} is obtained by identifying two  
module homomorphisms if their difference factors through a projective module.
Fortunately the two different meanings of `stable' fit together nicely; 
the stable module category is the homotopy category
associated to an underlying stable model category structure 
\cite[Sec.\ 2]{hovey-book}.  
Every finitely generated module is compact when considered as an object 
of the stable module category. Compare also 
Example \ref{single generator examples} (v).

\item {\bf Comodules over a Hopf-algebra.} Suppose $B$ is a  
commutative Hopf-algebra over a field. Hovey, Palmieri and  
Strickland introduce the category $\C(B)$ of chain complexes of  
injective $B$-comodules, with morphisms the chain homotopy classes  
of maps \cite[Sec.\ 9.5]{hps}. 
Compact generators are given by injective resolutions of simple comodules 
(whose isomorphism classes form a set). In \cite[Thm.\ 2.5.17]{hovey-book},  
Hovey shows that there is a cofibrantly generated model category  
structure on the category of {\em all} chain complexes of  
$B$-comodules whose homotopy category is the category $\C(B)$.
\end{enumerate}

\subsection{Quillen equivalences}\label{sec-Quillen-equiv} 
The most highly structured notion to express that two model
categories describe the same homotopy theory is that of a 
{\em Quillen equivalence}. 
An adjoint functor pair between model categories is a {\em Quillen pair} 
if the left adjoint $L$ preserves cofibrations
and trivial cofibrations. An equivalent condition is to demand that 
the right adjoint $R$ preserve fibrations and trivial fibrations. 
Under these conditions, the functors pass to an adjoint
functor pair on the homotopy categories, see \cite[I.4 Thm.\ 3]{Q},
\cite[Thm.\ 9.7 (i)]{DS} or \cite[1.3.10]{hovey-book}. A Quillen
functor pair is a {\em Quillen equivalence} if it induces an equivalence
on the homotopy categories. 
A Quillen pair is a Quillen equivalence if and only if the following 
criterion holds \cite[1.3.13]{hovey-book}: 
for every cofibrant object $A$ of the source category of $L$ 
and for every fibrant object $X$ of the source category of $R$,
a map $L(A)\to X$ is a weak equivalence if and only if its adjoint
$A \to R(X)$ is a weak equivalence. 

As pointed out in ~\cite[9.7 (ii)]{DS} and~\cite[I.4, Thm.\ 3]{Q}, 
in addition to inducing an equivalence of homotopy categories,
Quillen equivalences also preserve the homotopy theory
associated to a model category, that is, the higher order structure
such as mapping spaces, suspension and loop functors, and cofiber and fiber
sequences.
Note that the notions of compactness, generators, and stability 
are invariant under Quillen equivalences of model categories.
  
For convenience we restrict our attention to 
{\em simplicial} model categories (see \cite[II.2]{Q}). 
This is not a big loss of generality; it is shown in \cite{rss} that every 
cofibrantly generated, proper, stable model category is in fact 
Quillen equivalent to a simplicial model category. 
In \cite{dugger-simplicial}, Dugger obtains the same conclusion under somewhat
different hypotheses. In both cases the candidate is the category 
of simplicial objects over the given model category 
endowed with a suitable localization of the Reedy model structure.

\section{Classification theorems}
\label{sec-classification theorem}

\subsection{Monogenic stable model categories}

Several of the examples of stable model categories mentioned in 
\ref{stable examples} already come as categories of modules 
over suitable rings or ring spectra. This is no coincidence.
In fact, our first classification theorem  says that
every stable model category with a single compact
generator has the same homotopy theory as the modules over a
symmetric ring spectrum (see~\cite[5.4]{hss} 
for background on symmetric ring spectra). 
This is analogous to the classical fact~\cite[V1, p.\ 405]{gabriel} 
that module categories are characterized as those cocomplete
abelian categories which posses a single small projective  
generator; the classifying ring is obtained as the endomorphism ring 
of the generator.

In Definition \ref{def-endo ring spectrum} we associate to every object $P$ 
of a simplicial, cofibrantly generated, stable model
category $\C$ a symmetric {\em endomorphism ring spectrum} $\End(P)$.
The ring of homotopy groups $\pi_*\End(P)$ 
is isomorphic to the ring of graded self maps of $P$ in 
the homotopy category of $\C$, $[P,P]_{\ast}^{\HoC}$.

For the following theorem we have to make two technical assumptions. 
We need the notion of {\em cofibrantly generated
model categories} from~\cite{DHK} which is reviewed in some detail 
in~\cite[Sec.\ 2]{ss} and~\cite[Sec.\ 2.1]{hovey-book}. 
We also need properness (see \cite[Def.\ 1.2]{BF} or \cite[Def.\ 5.5.2]{hss}).
A model category is {\em left proper} if pushouts across  cofibrations
preserve weak equivalences.  A model category is {\em right proper}
if pullbacks over fibrations preserve weak equivalences.  A {\em proper}
model category is one which is both left and right proper.  

\begin{theorem} \label{thm-main-one} {\bf (Classification of  
monogenic stable model categories)} \\
Let $\C$ be a simplicial, cofibrantly generated,
proper, stable model category with a compact generator $P$.
Then there exists a chain of simplicial Quillen equivalences
between $\C$ and the model category of $\End(P)$-modules.
\[ \C \simeq_Q \mbox{\em{mod-}}\End(P) \]
\end{theorem}

This theorem is a special case of the more general
classification result Theorem  \ref{thm-main-gen}, which applies to 
stable model categories with a set of compact generators 
and which we prove in Section~\ref{symmetric objects}.
Furthermore if in the situation of Theorem \ref{thm-main-one}, 
$P$ is a compact object but not necessarily a generator of $\C$, 
then $\C$ still `contains' the homotopy theory of $\End(P)$-modules,
see Theorem \ref{main theorem spectral version} (ii) for the precise
statement.
In the Morita context (Theorem \ref{thm-Morita}) we also prove a partial
converse to Theorem \ref{thm-main-one}.

\subsection{Examples: stable model categories with a compact generator.}
\label{single generator examples}

\renewcommand{\labelenumi}{(\roman{enumi})}
\begin{enumerate}

\item {\bf Uniqueness results for stable homotopy theory.}
The classification theorem above yields a characterization 
of the model category of spectra: a simplicial, cofibrantly generated,
proper, stable model category is simplicially Quillen equivalent
to the category of symmetric spectra if and only if it has 
a compact generator $P$ for which the unit map of ring spectra
$S\to\End(P)$ is a stable equivalence. 
The paper \cite{ss-uniqueness} is devoted to other necessary and 
sufficient conditions for when a stable model category is 
Quillen equivalent to spectra -- some of them in terms of the
homotopy category of $\C$ and the natural action of the stable homotopy
groups of spheres. 
In \cite{sch-2}, this result is extended to a uniqueness theorem
showing that the 2-local stable homotopy category has only one
underlying model category up to Quillen equivalence. 
In both of these papers, we eliminate the technical conditions 
`cofibrantly generated' and `proper' by working with spectra in the sense of 
Bousfield and Friedlander \cite{BF},
as opposed to the Quillen equivalent symmetric spectra and 
`simplicial' by working with {\em framings} \cite[Chpt.~5]{hovey-book}.
In another direction, the uniqueness result is extended to include the
monoidal structure in \cite{ship-mon}.

\item {\bf Chain complexes and Eilenberg-Mac Lane spectra.} 
Let $A$ be a ring. Theorem~\ref{thm-chains and EM} 
shows that the model category of chain complexes of $A$-modules 
is Quillen  equivalent to the model category 
of modules over the symmetric Eilenberg-Mac Lane ring spectrum $HA$. 
This can be viewed as an instance of Theorem \ref{thm-main-one}: 
the free $A$-module of rank one, considered as a complex concentrated in 
dimension zero, is a compact generator for the unbounded derived category 
of $A$.
Since the homotopy groups of its endomorphism ring spectrum 
(as an object of the model category of chain complexes) are concentrated
in dimension zero, the endomorphism ring spectrum is stably equivalent to the 
Eilenberg-Mac Lane ring spectrum for $A$ 
(see Proposition \ref{prop-EMuniqueness}).  
This also shows that although the model category of chain complexes 
of $A$-modules is not simplicial it is Quillen equivalent 
to a simplicial model category.  So although our classification
theorems do not apply directly, they do apply indirectly.  

\item {\bf Smashing Bousfield localizations.} 
Let $E$ be a spectrum and consider the $E$-local model  
category structure on some model category of spectra 
(see e.g.\ \cite[VIII 1.1]{ekmm}).
This is another stable model category in which the localization of  
the sphere spectrum $L_E S^0$ is a generator. This localized sphere  
is compact if the localization is {\em smashing}, i.e., if a certain  
natural map $X\sm L_E S^0 \to L_E X$ is a stable equivalence for  
all $X$. So for a smashing localization the $E$-local model category  
of spectra is Quillen equivalent to modules over the ring spectrum  
$L_E S^0$ (which is the endomorphism ring spectrum of the localized  
sphere in the localized model structure).

\item {\bf $K(n)$-local spectra.} Even if a  
Bousfield localization is not smashing, Theorem \ref{thm-main-one}  
might be applicable. As an example we consider Bousfield  
localization with respect to the $n$-th Morava K-theory $K(n)$ at a  
fixed prime.
The localization of the sphere is still a generator,
but for $n>0$
it is not compact in the local category, see~\cite[3.5.2]{hps}.
However the localization of any finite type $n$ spectrum $F$ is a  
compact generator for the $K(n)$-local  
category~\cite[7.3]{hovey-strickland-kn}. Hence the $K(n)$-local  
model category is Quillen
equivalent to modules over the endomorphism ring $\End(L_{K(n)}F)$.

\item {\bf Frobenius rings.}
As in Example \ref{stable examples} (iv) we consider  
a Frobenius ring and assume that the stable module category has a  
compact generator. Then we are in the situation of Theorem  
\ref{thm-main-one}; however this example is completely algebraic,  
and there is no need to consider ring spectra to identify the stable  
module category as the derived category of a suitable `ring'. 
In fact Keller shows \cite[4.3]{keller-derivingDG} that in such a  
situation there exists a differential graded algebra (DGA) and  
an equivalence between the stable module category and the unbounded  
derived category of the DGA. 

A concrete example of this situation arises for group algebras of $p$-groups 
over a field $k$ of characteristic $p$. 
In this case the trivial module is the only simple module, and it is a
compact generator of the stable module category. 
More generally a result of Benson \cite[Thm.\ 1.1]{benson-principal_block}
says that the trivial module generates the stable module category
of the principal block of a group algebra $kG$ if and only if 
the centralizer of every element of order $p$ is  $p$-nilpotent. 
So in this situation Keller's theorem applies and identifies the stable module
category as the unbounded derived category of a certain DGA.
The homology groups of this DGA are isomorphic (by construction) to the
ring of graded self maps of the trivial module in the stable module category, 
which is just the Tate-cohomology ring $\widehat{\mbox{H}}^*(G; k)$.

\item \label{ex-t-alg} {\bf Stable homotopy of algebraic theories.}  
Another motivation for this paper and an early instance of Theorem  
\ref{thm-main-one} came from the stabilization of the model category  
of algebras over an algebraic theory \cite{sch-theories}.
For every pointed algebraic theory $T$, the category of simplicial
$T$-algebras is a simplicial model category so that one has a 
category $\Spc(T)$ of (Bousfield-Friedlander type) spectra 
of $T$-algebras, a cofibrantly generated,  simplicial stable
model category \cite[4.3]{sch-theories}. 
The free $T$-algebra on one generator has an endomorphism ring
spectrum which is constructed as a Gamma-ring in \cite[4.5]{sch-theories}
and denoted $T^s$. Then \cite [Thm.\ 4.4]{sch-theories} provides 
a Quillen equivalence between the categories of connective spectra 
of $T$-algebras and the category of $T^s$-modules 
(the connectivity condition could be removed by working 
with symmetric spectra instead of $\Gamma$-spaces). 
This fits with Theorem  \ref{thm-main-one} because
the suspension spectrum of the free $T$-algebra on one generator is 
a compact generator for the category $\Spc(T)$. 
See \cite[Sec.\ 7]{sch-theories} for a list of ring spectra that 
arise from algebraic  theories in this fashion.

\end{enumerate}

\begin{remark} \label{rem-K(n)}
{\em The notion of a compact generator and the  
{\em homotopy groups} of the endomorphism ring spectrum only depend 
on the homotopy category, and so they are invariant under  
equivalences of triangulated categories. However, the {\em homotopy type}  
of the endomorphism ring spectrum depends on the model category  
structure. The following example illustrates this point.  Consider  
the $n$-th Morava K-theory spectrum $K(n)$ for a fixed prime and  
some number $n> 0$. This spectrum admits the structure of an  
A$_{\infty}$-ring spectrum \cite{robinson-K(n)}.  
Hence it also has a model as an $S$-algebra or a symmetric ring spectrum 
and the category of its module spectra is a stable model category. 
The ring of homotopy groups of $K(n)$ is the graded field 
${\mathbb F}_p[v_n,v_n^{-1}]$ with $v_n$  
of dimension $2p^n-2$. Hence the homotopy group functor establishes  
an equivalence between the homotopy category of $K(n)$-module spectra  
and the category of graded ${\mathbb F}_p[v_n,v_n^{-1}]$-modules.

Similarly the homology functor establishes an equivalence 
between the derived  category of differential graded modules over the 
graded field ${\mathbb F}_p[v_n,v^{-1}_n]$ and the category of
graded  ${\mathbb F}_p[v_n,v^{-1}_n]$-modules.  
So the two stable model categories of $K(n)$-module spectra and 
DG-modules over ${\mathbb F}_p[v_n,v^{-1}_n]$
have equivalent triangulated homotopy categories 
(including the action of the stable homotopy groups of spheres 
--- all elements in positive dimension act trivially in both cases). 
But the endomorphism ring spectra of the respective free rank one modules 
are the Morava $K$-theory ring spectrum on the one side 
and the Eilenberg-Mac Lane ring spectrum for ${\mathbb F}_p[v_n,v^{-1}_n]$
on the other side, which are not stably equivalent. 
Similarly the two model categories are not Quillen
equivalent since for DG-modules all function spaces are products 
of Eilenberg-Mac Lane spaces, but for $K(n)$-modules they are not.
}
\end{remark}

\subsection{Multiple generators}

There is a generalization of Theorem \ref{thm-main-one} to the
case of a stable model category with a set of compact generators
(as opposed to a single compact generator).

Let us recall the algebraic precursors of this result: a {\em ringoid}
is a category whose hom-sets are abelian groups with bilinear composition. 
Ringoids are sometimes called {\em pre-additive categories} 
or {\em rings with several objects}.
Indeed a ring in the traditional sense is the same as a ringoid with 
one object. 
A (right) {\em module} over a ringoid is defined to
be a contravariant additive functor to the category of abelian groups. 
These more general module categories have been identified 
as the cocomplete abelian categories 
which have a set of small projective generators \cite[5.3H, p. 120]{freyd}. 
An analogous theory for derived
categories of DG categories has been developed by Keller 
\cite{keller-derivingDG}.
 
Our result is very much in the spirit of Freyd's or Keller's, 
with spectra substituting for abelian groups 
or chain complexes.
A symmetric ring spectrum can be viewed as a category with
one object which is enriched over symmetric spectra; 
the module category then becomes the category of enriched (spectral) functors 
to symmetric spectra. So we now look at `ring spectra with several
objects' which we call {\em spectral categories}. 
This is analogous to pre-additive, differential graded
or simplicial categories which are enriched over abelian groups, 
chain complexes or simplicial sets respectively.

\begin{definition} \label{def-spectral cat} {\em A {\em spectral category}
is a category $\Oc$ which is enriched over the category $\spec$ 
of symmetric spectra (with respect to smash product, 
i.e., the monoidal closed structure of \cite[2.2.10]{hss}). 
In other words, for every pair of objects $o,o'$ in $\Oc$ there is a 
morphism symmetric spectrum $\Oc(o,o')$, 
for every object $o$ of $\Oc$ there is a map 
from the sphere spectrum $S$ to $\Oc(o,o)$ (the `identity element' of $o$),
and for each triple of objects there is an associative and unital 
composition map of symmetric spectra
\[ \Oc(o',o'') \ \sm \ \Oc(o, o') \ \varrow{1cm} \ \Oc(o,o'') \ . \]
An $\Oc$-module $M$ is a contravariant spectral functor to the category 
$\spec$  of symmetric spectra,
i.e., a symmetric spectrum $M(o)$ for each object of $\Oc$ 
together with coherently associative and unital maps of symmetric spectra
\[ M(o) \ \sm \ \Oc(o',o) \ \varrow{1cm} \  M(o') \]
for pairs of objects $o,o'$ in $\Oc$.
A morphism of $\Oc$-modules $M\to N$ consists of maps of symmetric
spectra $M(o) \to N(o)$ strictly compatible with the action of $\Oc$. 
We denote the category of $\Oc$-modules by mod-$\Oc$. 
The {\em free} (or `representable') module $F_o$ is given by 
$F_o(o')=\Oc(o',o)$.}
\end{definition}

\begin{remark}{\em In Definition \ref{def-spectral cat} 
we are simply spelling out what it means to do enriched
category theory over the symmetric monoidal closed category $\spec$ 
of symmetric spectra with respect to the smash product 
and the internal homomorphism spectra. Kelly's book \cite{kelly} is
an introduction to enriched category theory in general; 
the spectral categories, modules over these (spectral functors) and 
morphisms of modules 
as defined above are the $\spec$-categories, $\spec$-functors and
$\spec$-natural transformations in the sense of \cite[1.2]{kelly}. 
So the precise meaning of the coherence
and compatibility conditions in Definition \ref{def-spectral cat} can be found 
in \cite[1.2]{kelly}.}
\end{remark}

We show in Theorem \ref{O-modules} that for any spectral category $\Oc$
the category of $\Oc$-modules is a model category with objectwise
stable equivalences as the weak equivalences.  There we also show
that the set of free modules $\{F_o\}_{o\in \Oc}$ is a set
of compact generators.  If $\Oc$ has a single object $o$,
then the $\Oc$-modules are precisely the modules over the symmetric 
ring spectrum $\Oc(o,o)$, and the model category structure is  
the one defined in \cite[5.4.2]{hss}.

In Definition \ref{def-endo ring spectrum} we associate to every set $\Pc$ 
of objects of a simplicial, cofibrantly generated stable model
category $\C$ a spectral {\em endomorphism category} $\EP$ whose objects
are the members of the set $\Pc$ and such that
there is a natural, associative and unital 
isomorphism 
\[ \pi_* \, \EP(P,P') \ \iso \ [P,P']_*^{\HoC} \ . \]
For a set with a single element this reduces to the notion of the 
endomorphism ring spectrum.

\begin{theorem} \label{thm-main-gen} 
{\bf (Classification of stable model categories)}
Let $\C$ be a simplicial, cofibrantly generated, proper, 
stable model category with a set $\Pc$ of compact generators.
Then there exists a chain of simplicial Quillen equivalences 
between $\C$ and the model category of $\EP$-modules.
\[ \C \simeq_Q \mbox{{\em mod-}}\EP \]
\end{theorem}

There is an even more general version of Theorem
\ref{thm-main-gen} which also provides information if the set $\Pc$
does not generate the whole  homotopy category, 
see Theorem \ref{main theorem spectral version} (ii). 
This variant implies that for any set $\Pc$ of compact objects in a proper,
cofibrantly generated, simplicial, stable model category 
the homotopy category of $\EP$-modules is triangulated equivalent 
to the localizing subcategory of $\HoC$ generated by  the set $\Pc$.

The  proof of Theorem \ref{thm-main-gen} breaks up into two parts.  
In order to mimic the classical proof for abelian categories  
we must consider a situation where the hom functor $\Hom_{\C}(P,-)$ 
takes values in the category of modules over a suitable 
endomorphism ring spectrum of $P$.
In Section \ref{symmetric objects} we show how this can be arranged, 
given the technical conditions that $\C$ is cofibrantly generated, 
proper and simplicial. 
We introduce the category $\specC$ of symmetric spectra over $\C$ 
and show in Theorem \ref{thm-specC} that it is a stable model category 
with composable and homotopically well-behaved function symmetric spectra 
which is Quillen equivalent to the original stable model category $\C$. 

In Theorem \ref{main theorem spectral version} we prove 
Theorem \ref{thm-main-gen} under the assumption that 
$\C$ is a {\em spectral model category} (Definition 
\ref{def-spectral model cat}), i.e., a model category 
with composable and homotopically well-behaved function symmetric spectra.
Since the model category $\specC$ is spectral and Quillen equivalent to $\C$ 
(given the technical assumptions of Theorem \ref{thm-main-gen}),
this will complete the classification results.  

\subsection{Examples: stable model categories with a set of generators.}
\label{multiple generator examples} 

\renewcommand{\labelenumi}{(\roman{enumi})}
\begin{enumerate}
\item {\bf Equivariant stable homotopy.}
Let $G$ be a compact Lie group. As mentioned in 
Example~\ref{stable examples}(iii)
there are several versions of model categories of $G$-equivariant
spectra. In~\cite{mm}, $G$-equivariant orthogonal spectra are shown
to form a cofibrantly generated, topological (hence simplicial),
proper model category (the monoidal structure 
plays no role for our present considerations).
For every closed subgroup $H$ of $G$ the equivariant suspension spectrum 
of the homogeneous space $G/H_+$ is compact 
and the set $\G$ of these spectra for all closed subgroups $H$ 
generates the $G$-equivariant stable homotopy category, see \cite[9.4]{hps}.
      
Recall from \cite[V \S 9]{lms} that a {\em Mackey functor} 
is a module over the stable homotopy  orbit category, 
i.e., an additive functor from the homotopy orbit category $\pi_0\,\E(\G)$ 
to the category of abelian groups; by \cite[V Prop.\ 9.9]{lms} 
this agrees with the original algebraic definition of Dress \cite{dress} 
in the case of finite groups. The spectral endomorphism category $\E(\G)$ is 
a spectrum valued lift of the stable homotopy orbit category. 
Theorem \ref{thm-main-gen} shows  that the category of $G$-equivariant 
spectra is Quillen equivalent to a category of {\em topological
Mackey functors}, i.e., the category of modules over the stable orbit 
category $\E(\G)$. Note that the homotopy type of each morphism spectrum 
of $\E(\G)$ depends on the universe $\U$.  

After rationalization the Mackey functor analogy becomes even more concrete:
we will see in Example \ref{ex-rat equivariant via tilting}
that the model category of rational $G$-equivariant spectra is in fact Quillen
equivalent to the model category of chain complexes of 
rational Mackey functors for $G$ finite.
For certain non-finite compact Lie groups, our approach via `topological
Mackey functors' is used in~\cite{shipley-S1} and~\cite{greenlees-shipley} 
as an intermediate step in forming algebraic models 
for rational $G$-equivariant spectra.

\item {\bf Motivic stable homotopy of schemes.} 
In \cite{morel-voevodsky,voevodsky-icm} Morel and Voevodsky introduce
the ${\mathbb A}^1$-local model category structure for schemes over a base. 
The objects  of their category are sheaves of sets in the Nisnevich topology 
on smooth schemes of finite type over a fixed base scheme. 
The weak equivalences are the ${\mathbb A}^1$-local equivalences 
-- roughly speaking they are generated by the projection maps 
$X\times {\mathbb A}^1 \to X$ for smooth schemes $X$,
where ${\mathbb A}^1$ denotes the affine line.

Voevodsky \cite[Sec.\ 5]{voevodsky-icm} introduces an associated 
stable homotopy category by inverting smashing with the `Tate-sphere' 
$T = {\mathbb A}^1/({\mathbb A}^1- 0)$. The punch-line is that  theories
like algebraic $K$-theory or motivic cohomology are represented by objects in
this stable homotopy category \cite[Sec.\ 6]{voevodsky-icm},
 at least when the base scheme is the spectrum of a field.
      
In \cite{jardine-A^1stable}, Jardine provides the details of the  construction
of model categories of $T$-spectra over the spectrum of a field $k$. 
He constructs two Quillen equivalent proper,  simplicial model categories 
of Bousfield-Friedlander type and symmetric ${\mathbb A}^1$-local $T$-spectra
\cite[2.11, 4.18]{jardine-A^1stable}. Since $T$ is weakly equivalent 
to a suspension (of the  multiplicative group scheme), 
this in particular yields a stable model category. 
A set of  compact generators for this homotopy category is given 
by the $T$-suspension spectra $\Sigma^{\infty}_T(\mbox{Spec}R)_+$ 
when $R$ runs over smooth $k$-algebras of finite type. 
So if $k$ is countable then this is a countable set of compact generators, 
compare \cite[Prop.\ 5.5]{voevodsky-icm}.

\item {\bf Algebraic examples.} 
Again the classification theorem \ref{thm-main-gen} has an  
algebraic analogue and precursor, namely Keller's theory of derived  
equivalences of DG categories \cite{keller-derivingDG}. The bottom  
line is that if an example of a stable model category is algebraic  
(such as derived or stable module categories in Examples 
\ref{abelian examples}), then it is not  
necessary to consider spectra and modules over spectral categories,  
but one can work with chain complexes and differential graded  
categories instead. As an example, Theorem 4.3 of  
\cite{keller-derivingDG} identifies the stable module category of a  
Frobenius ring with the unbounded derived category of a certain  
differential graded category.
\end{enumerate}

\subsection{Prerequisites on spectral model categories.}
\label{sub-prerequisites}

A spectral model category is analogous to a simplicial model category,
\cite[II.2]{Q}, but with the category of simplicial sets replaced 
by symmetric spectra.  Roughly speaking, a spectral
model category is a pointed model category which is compatibly 
enriched over the stable model category of symmetric spectra. 
The compatibility is expressed by the axiom (SP) below
which takes the place of \cite[II.2 SM7]{Q}. 
For the precise meaning of `tensors' and `cotensors'
over symmetric spectra see e.g.\ \cite[3.7]{kelly}. 
A spectral model category is the  same as a `$\spec$-model category' 
in the sense of \cite[Def.\ 4.2.18]{hovey-book}, where the category 
of symmetric spectra is endowed with the stable model structure of 
\cite[3.4.4]{hss}. Condition two of \cite[4.2.18]{hovey-book} is
automatic since the unit $S$ for the smash product of symmetric spectra
is cofibrant.
Examples of spectral model categories are module categories over
a symmetric ring spectrum, 
module categories over a spectral category (Theorem \ref{O-modules}) 
and the category of symmetric spectra over a suitable simplicial model category
(Theorem \ref{thm-specC}).

\begin{definition} \label{def-spectral model cat}
{\em A {\em spectral model category} is a model category $\C$ which is 
tensored, cotensored and enriched (denoted $\Hom_{\C}$) over the category 
of symmetric spectra with the closed monoidal
structure of \cite[2.2.10]{hss} such that the following compatibility 
axiom (SP) holds:\\
(SP) For every cofibration $A\to B$ and every fibration $X\to Y$ in $\C$ 
the induced map
\[ \Hom_{\C}(B,X) \ \varrow{1cm} \ \Hom_{\C}(A,X) \times_{\Hom_{\C}(A,Y)}
\Hom_{\C}(B,Y) \]
is a stable fibration of symmetric spectra. If in addition one of the maps 
$A\to B$ or $X\to Y$ is a weak equivalence, then the resulting map 
of symmetric spectra is also a stable equivalence.
We use the notation $K\sm X$ and $X^K$ to denote the tensors and cotensors
 for $X\in \C$ and $K$ a symmetric spectrum.}
\end{definition}

In analogy with \cite[II.2 Prop.\ 3]{Q} the compatibility axiom (SP) 
in Definition \ref{def-spectral model cat} of a spectral model category 
can be cast into two adjoint forms, one of which will be of use for us.  
Given a categorical enrichment of a model category $\C$ over the category 
of symmetric spectra, then axiom (SP) is equivalent to (SPb) below.  
The equivalence of conditions (SP) and (SPb) is a consequence of the 
adjointness properties of the tensor and cotensor functors, 
see \cite[Lemma 4.2.2]{hovey-book} for the details.

(SPb) For every  
cofibration $A\to B$ in $\C$ and every stable cofibration $K\to L$  
of symmetric spectra, the canonical map ({\em pushout product map})
\[ L \sm A \cup_{K \sm A}  K \sm B \ \varrow{1cm} \ L \sm B \]
is a cofibration; the pushout product map is a weak equivalence if  
in addition $A\to B$ is a weak equivalence in $\C$ or $K\to L$ is a  
stable equivalence of symmetric spectra.
  
In a spectral model category the levels of the symmetric function spectra 
$\Hom_{\C}(X,Y)$ can be rewritten as follows. 
The adjunctions give an isomorphism of simplicial sets
\[  \Hom_{\C}(X,Y)_n \ \iso \ \map_{\Spc}(F_nS^0,\Hom_{\C}(X,Y)) \ \iso \
\map_{\C}(F_nS^0 \sm X,Y) \ \iso \ \map_{\C}(X,Y^{F_nS^0}) \]
where $F_nS^0$ is the free symmetric spectrum generated at level $n$ 
by the 0-sphere  (see \cite[2.2.5]{hss} or Definition \ref{free}).

\begin{lemma} \label{lemma-spectral implies stable}
A spectral model category is in particular a {\em simplicial}
and {\em stable}  model category. For $X$ a cofibrant and $Y$ a 
fibrant object of a spectral model category $\C$  there is a natural
isomorphism of graded abelian groups
$\pi^s_* \Hom_{\C}(X,Y) \iso [X,Y]^{\HoC}_*$.
\end{lemma}
\begin{proof}
The tensor and cotensor of an object of $\C$ 
with a pointed simplicial set $K$ 
is defined by applying the tensor and cotensor with the symmetric 
suspension spectrum $\Sigma^{\infty} K$. The homomorphism simplicial set 
between two objects of $\C$ is the 0-th level of the homomorphism symmetric
spectrum. The necessary adjunction formulas and the compatibility axiom 
\cite[II.2 SM7]{Q}
hold because the suspension spectrum functor 
$\Sigma^{\infty}\colon\sset\to \Spc$ from the category of pointed simplicial
sets to symmetric spectra is the left adjoint
of a Quillen adjoint functor pair and preserves the smash product 
(i.e., it is strong symmetric monoidal).  
In order to see that $\C$ is stable we recall that the shift functor in 
the homotopy category of $\C$ is the suspension functor. 
For cofibrant objects suspension is represented on the
model category level by the smash product with the one-dimensional 
sphere spectrum $\Sigma^{\infty} S^1$.  
This sphere spectrum is invertible, up to stable 
equivalence of symmetric spectra, with inverse the
(-1)-dimensional sphere spectrum (modeled as a symmetric spectrum by
$F_1 S^0$).  
Since the action of symmetric spectra 
on $\C$ is associative up to coherent isomorphism this implies 
that suspension is a self-equivalence of the homotopy
category of $\C$.  This in turn implies that the right adjoint 
loop functor has to be an inverse equivalence. 
If $X$ is cofibrant and $Y$ fibrant in $\C$, then by the 
compatibility axiom (SP) the symmetric spectrum $\Hom_{\C}(X,Y)$ 
is stably fibrant, i.e., an $\Omega$-spectrum. 
So for $n\geq  0$ we have isomorphisms
\[\pi_n \, \Hom_{\C}(X,Y) \ \iso \ \pi_n \map_{\C}(X,Y) \ \iso \ \pi_0 
\map_{\C}(X,\Omega^nY) \ \iso \  [X,\Omega^ nY]^{\HoC} \ ; \]
and for $n \leq 0$ we have the isomorphisms
\[ \pi_n \Hom_{\C}(X,Y) \ \iso  \ \pi_0 \Hom_{\C}(X,Y)_n \ \iso \ 
\pi_0 \map_{\C}(F_nS^0 \sm X,Y)  \ \iso \ [\Omega^nX,Y]^{\HoC} \ . \]
\end{proof}

\subsection{Symmetric spectra over a category}
\label{symmetric objects}

Throughout this section we assume that $\C$ is a cocomplete  
category which is tensored and cotensored over the category $\sset$ of
pointed simplicial sets, with this action denoted by $\tensor$ and  
morphism simplicial sets denoted $\map_{\C}$.
We let $S^1= \Delta[1]/\partial\Delta[1]$ be our model for the  
simplicial circle and we set $S^n = (S^1)^{\wedge n}$ for $n>1$; the  
symmetric group on $n$ letters acts on $S^n$ by permuting the  
coordinates.

\begin{definition} \label{def-specC} 
{\em Let $\C$ be a category which is tensored  
over the category of pointed simplicial sets. 
A {\em symmetric sequence over $\C$} is a sequence of objects 
$X=\{X_n\}_{n\geq 0}$ in $\C$ together with a left action of 
the symmetric group $\Sigma_n$ on  $X_n$ for all $n\geq 0$. 
A {\em symmetric spectrum over $\C$} is a symmetric  sequence 
in $\C$ with coherently associative  
$\Sigma_p\times\Sigma_q$-equivariant morphisms
\[ S^p \, \tensor \, X_q \ \varrow{1cm} \ X_{p+q}  \]
(for all $p,q\geq 0$). A morphism of symmetric sequences or 
symmetric spectra $X\to Y$ consists of a sequence of $\Sigma_n$-equivariant 
morphisms $X_n\to Y_n$ which commute with the structure maps.
We denote the category of symmetric sequences by $\C^{\Sigma}$ 
and the category of symmetric spectra by $\specC$.}
\end{definition}

Since $\C$ is simplicial, the action of $\sset$ on $\C$ 
extends to an action of $\sset^{\Sigma}$, the category of symmetric 
sequences over $\sset$, on $\C^{\Sigma}$:

\begin{definition}\label{def-tensor-sset}
{\em Given $X$ a symmetric sequence over $\C$, and $K$ a  
symmetric sequence over
$\sset$, we define their {\em tensor product}, $K\tensor X$,
by the formula
\[ (K\tensor X)_{n}=\bigvee _{p+q=n} \Sigma _{n}^{+}
\tensor_{\Sigma_{p}\times \Sigma _{q}} (K_{p}\tensor X_{q})\ . \]}
\end{definition}

The symmetric sequence $S=(S^0, S^1, \cdots, S^n, \cdots)$ of
simplicial sets is a commutative monoid in
the symmetric monoidal category $(\sset^{\Sigma},\tensor)$.
The unit map here is the identity in the first spot and the base point
elsewhere, $\eta \colon (S^0, *, *, \cdots)=u \to S$. In this  
language a symmetric spectrum over $\C$ can be redefined as a left  
$S$-module in the category of symmetric sequences over~$\C$.

\begin{definition}
{\em Let $X$ be an object in $\specC$ and $K$ a symmetric spectrum in $\spec$.
Define their  {\em smash product}, $K \wedge X$, as the
coequalizer of the two maps
\[ K \tensor S \tensor X \ \parallelarrows{0.5cm} \ K \tensor X \ \]
induced by the action of $S$ on $X$ and $K$ respectively. 
So $\specC$ is tensored over the symmetric monoidal category 
of symmetric spectra.
Dually, we define a symmetric spectrum valued morphism object
$\Hom_{\specC}(X,Y) \in \spec$ for $X,Y \in \specC$.  As a preliminary
step, define a shifting down functor, $\sh_n\colon \specC \xrightarrow{} 
\specC$, by $(\sh_nX)_m= X_{n+m}$ where $\Sigma_m$ acts via the inclusion 
into $\Sigma_{n+m}$. 
Note there is a leftover action of $\Sigma_n$ on $\sh_n X$.   
Define $\Hom_{\Sigma}(X,Y) \in \sset^{\Sigma}$ for $X,Y$ objects 
in $\C^{\Sigma}$ by $\Hom_{\Sigma}(X,Y)_n=\map(X, \sh_n Y)$, the
simplicial mapping space given by the simplicial structure on  $\C^{\Sigma}$
with the $\Sigma_n$ action given by the leftover action of $\Sigma_n$ on
$\sh_n Y$ as mentioned above.
Then $\Hom_{\specC}(X,Y) \in \spec$ is the equalizer of the two maps 
$\Hom_{\Sigma}(X,Y) \to \Hom_{\Sigma}(S\tensor X, Y)$.
}
\end{definition}

Using this spectrum valued hom functor, for $K$ in $\spec$ and $Y$ in $\specC$,
we define $Y^K \in \specC$ as the adjoint of the functor $K\tensor -$. 
That is, for any $X \in \specC$ define $Y^K$ such that
\begin{equation} \label{spectral adjunction}
\Hom_{\specC}(K\tensor X,Y) \iso \Hom_{\specC}(X,Y^K)
\iso \Hom_{\spec}(K,\Hom_{\specC}(X,Y)) \ .
\end{equation}

\begin{definition}\label{free}
{\em The $n$th evaluation functor $\Ev_n: \specC \to \C$ 
is given by $\Ev_n(X)=X_n$, ignoring the action of the symmetric group. 
The functor $\Ev_n$ has a left adjoint $F_n\colon \C \to \specC$ 
which has the form 
$(F_nX)_m=  \Sigma_{m}^+\tensor_{\Sigma_{m-n}} S^{m-n}\tensor X$ 
where $S^n=*$ for $n < 0$.
We use $\Sigma^{\infty}$ as another name for $F_0$ and call it the 
{\em suspension spectrum}.}
\end{definition}

\subsection{The level model structure on $\specC$}\label{level model}

There are two model category structures on symmetric spectra over $\C$  
which we consider; the level model category which we discuss in this
section, and the stable model category 
(see Section \ref{sub-stable structure}).
The level model category is a stepping stone 
for defining the stable model category, but it also allows us to define
endomorphism ring spectra (Definition \ref{def-endo ring spectrum}).

\begin{definition} {\em Let $f:X \to Y$ be a map in $\specC$.
The map $f$ is a {\em level equivalence} if each $f_n:X_n \to Y_n$ 
is a weak equivalence in $\C$, ignoring the $\Sigma_n$ action.
It is a {\em level fibration} if each $f_n$ is a fibration in $\C$. 
It is a  {\em cofibration} if it has the
left lifting property with respect to all level trivial fibrations.} 
\end{definition}

\begin{proposition}\label{prop-level-model}
For any simplicial, cofibrantly generated model category $\C$, 
$\specC$ with the level equivalences, level fibrations, and cofibrations
described above forms a cofibrantly generated model category referred 
to as the {\em level model category}, and denoted by  $\specC\rscript{\em lv}$.
Furthermore the following level analogue of the spectral axiom 
{\em (SP)} holds: \smallskip \\
{\em (SP$\rscript{lv}$)} for  every cofibration 
$A\to B$ and every level fibration $X\to Y$ in $\specC$ the induced map
\[ \Hom_{\specC}(B,X) \ \varrow{1cm} \ \Hom_{\specC}(A,X) 
\times_{\Hom_{\specC}(A,Y)} \Hom_{\specC}(B,Y) \]
is a level fibration of symmetric spectra. If in addition one of the maps 
$A\to B$ or $X\to Y$ is a level equivalence, then the resulting map 
of symmetric spectra is also a level equivalence.
\end{proposition}

There are various theorems in the model category literature which are useful 
in establishing model category structures. These theorems separate 
formal considerations that tend to show up routinely from 
the properties which require special arguments in each specific case. 
Since we will construct model category structures several times in this paper,
we recall  one such result that we apply in our cases. 
We work with the concept of cofibrantly generated model categories,
introduced by Dwyer, Hirschhorn and Kan \cite{DHK}; 
see also Section 2.1 of Hovey's book \cite{hovey-book} 
for a detailed treatment of this concept.
  
We use the same terminology as \cite[Sec.\ 2.1]{hovey-book}. 
Let $I$ be a set of maps in a  category. A map is a 
{\em relative $I$-cell complex} if it is a (possibly transfinite) 
composition of cobase changes of maps in $I$. An {\em $I$-injective} map 
is a map with the right lifting property with respect to  every map in $I$.
An {\em $I$-cofibration} is a map with the left lifting property 
with respect to $I$-injective maps.
For the definition of smallness relative to a set of maps 
see~\cite[2.1.3]{hovey-book}.

One of the main properties of cofibrantly generated model categories is that 
they admit an abstract version of Quillen's small object argument 
\cite[II 3.4]{Q}.

\begin{lemma} \label{transfinite small object}
{\em \cite{DHK}, \cite[2.1.14, 2.1.15]{hovey-book}}
Let $\C$ be a cocomplete category and $I$ a set of
maps in $\C$ whose domains are small relative to the relative 
$I$-cell complexes. Then
\begin{itemize}
\item there is a functorial factorization of any map $f$ in $\C$ 
as $f = qi$ with $q$ an $I$-injective map and $i$ a relative 
$I$-cell complex, and thus
\item  every $I$-cofibration is a retract of a relative $I$-cell complex.
\end{itemize}
\end{lemma}

\begin{theorem} \label{thm-lifting model structure}
{\em \cite{DHK}, \cite[Thm.\ 2.1.19]{hovey-book}} 
Let $\C$ be a complete and cocomplete category   and $I$ and
$J$ two sets of maps of $\C$ such that the domains of the maps in $I$ 
and $J$ are small  with respect to the relative $I$-cell complexes 
and the relative $J$-cell complexes respectively. 
Suppose also that a subcategory of $\C$ is specified whose morphisms 
are called `weak equivalences'.
 
Then there is a cofibrantly generated model structure on $\C$ with the 
given class of weak equivalences, with $I$ a set of 
generating cofibrations, and with $J$ a set of 
generating trivial  cofibrations if the following conditions hold:
\renewcommand{\labelenumi}{(\arabic{enumi})}
\begin{enumerate}
\item if $f$ and $g$ are composable morphisms such that two of the 
three maps $f$, $g$ and $gf$ are weak equivalences, 
then the third is also a weak equivalence.
\item every relative $J$-cell complex is an $I$-cofibration 
and a weak equivalence.
\item  the $I$-injectives are precisely the maps which are 
both $J$-injective and weak  equivalences.
\end{enumerate}
\end{theorem}

\begin{proof}[Proof of Proposition \ref{prop-level-model}]
Let $I_{\C}$ and $J_{\C}$ be sets of generators for 
the cofibrations and trivial cofibrations of $\C$. 
We define sets of generators for the level model category by 
$FI_{\C} = \{F_nI_{\C}\}_{n\geq 0}$ and $FJ_{\C} = \{F_nJ_{\C}\}_{n\geq 0}$, 
i.e., $F_n$ applied to the generators of $\C$ for each $n$. 
Then the $FI_{\C}$-injectives are precisely the levelwise
trivial fibrations and the $FJ_{\C}$-injectives are precisely the 
level fibrations.
We claim that every relative $FI_{\C}$-cell complex is levelwise 
a cofibration in $\C$, and similarly every relative $FJ_{\C}$-cell complex is 
levelwise a trivial cofibration in $\C$.
We show this for relative $FJ_{\C}$-cell complexes; the argument for 
$FI_{\C}$ is the same. 
Since level evaluation preserves colimits it suffices to check the claim 
for the generating cofibrations, $F_n A \to F_nB$ for $A \to B \in J_{\C}$. 
But the $m$th level of this map is a coproduct of $m!/(m - n)!$
copies of the map $S^{m-n} \sm A \to  S^{m-n} \sm B$. 
By the simplicial compatibility axiom \cite[II.2 SM7]{Q}, smashing
with a simplicial sphere preserves trivial cofibrations, so we are done. 
 
Now we apply Theorem \ref{thm-lifting model structure} 
to the sets $FI_{\C}$ and $FJ_{\C}$ with the level equivalences as weak
equivalences. Checking that the maps in $FI_{\C}$ are small with respect to 
relative $FI_{\C}$-cell complexes comes down (by adjointness) to checking 
that the domains of the maps in $I_{\C}$ are small with respect to
the levels of relative $FI_{\C}$-cell complexes, i.e. the cofibrations
in $\C$.  By \cite[14.2.14]{hirschhorn-book}, since the domains in $I_{\C}$ are
small with respect to the relative $I_{\C}$-cell complexes they are also
small with respect to all cofibrations.
The argument for $FJ_{\C}$ is the same. 
Conditions (1) and (3) of Theorem \ref{thm-lifting model structure} 
hold and condition (2) follows from the above claim. 
So we indeed have a cofibrantly generated level model structure.
  
To prove the property (SP$\rscript{lv}$) it suffices 
to check its adjoint pushout product form, i.e., the level analogue of
condition (SPb) of Section \ref{sub-prerequisites}; 
it is enough to show that the  pushout product of two generating cofibrations 
is a cofibration, and similarly when one of the maps is a trivial cofibration
(see \cite[2.3 (1)]{ss} or \cite[Cor.\ 4.2.5]{hovey-book}). 
So let $i\in I_{\sset}$ and let $j \in  I_{\C}$. 
Then the product $F_n i\sm  F_m j$ is isomorphic to $F_{n+m}(i\sm j)$ and the 
result follows since the free functors preserve cofibrations and trivial 
cofibrations. 
\end{proof}

We can now introduce endomorphism ring spectra and endomorphism categories.

\begin{definition} \label{def-endo ring spectrum}
{\em Let $\C$ be a simplicial and cofibrantly generated model category 
and $\Pc$ a set of cofibrant objects. For every object $P\in \Pc$
let $\Sigf P$ be a fibrant replacement of the symmetric suspension
spectrum of $P$ in the level model structure on $\specC$ of 
Proposition \ref{prop-level-model}. We define the 
{\em endomorphism category} $\EP$  as the full spectral subcategory
of $\specC$ with objects $\Sigf P$ for $P\in \Pc$. 
To simplify notation, we usually denote 
objects of $\EP$ by $P$ instead of $\Sigf P$.
If $\Pc$ has a single object $P$ we also refer to the symmetric ring spectrum
$\EP(P,P)=\Hom_{\specC}(\Sigf P,\Sigf P)$ as the 
{\em endomorphism ring spectrum} of the object $P$.}
\end{definition}

Up to stable equivalence, the definition of the endomorphism category 
does not depend on the choices of fibrant replacements. 

\begin{lemma}
Let $\C$ be a simplicial and cofibrantly generated model 
category and $\Pc$ a set of cofibrant objects.
Suppose $\{\Sigf P\}_{P\in\Pc}$ and  $\{\widebar{\Sigf P}\}_{P\in\Pc}$ 
are two sets of level fibrant replacements 
of the symmetric suspension spectra. Then the two full spectral subcategories 
of $\specC$ with objects  $\{\Sigf P\}_{P\in\Pc}$ and  
$\{\widebar{\Sigf P}\}_{P\in\Pc}$ respectively are stably equivalent.
\end{lemma}
\begin{proof}
The proof uses the notion of {\em quasi-equivalence},
see Definition \ref{def-quasi-equivalence}. For every $P\in \Pc$ we choose
a level equivalence $\phi_P: \widebar{\Sigf P}\to \Sigf P$. 
We define a $\EP$-$\widebar{\EP}$-bimodule $M$ by the rule
\[ M(P,P') \ = \ \Hom_{\specC}(\widebar{\Sigf P},\Sigf P) \ . \]
Because of the property (SP$\rscript{lv}$) of the 
homomorphism spectra in $\specC$ the bimodule $M$ is a quasi-equivalence
with respect to the maps $\phi_P$, and the result follows 
from Lemma \ref{lemma-balanced}.
\end{proof}

\subsection{The stable model structure on $\specC$}
\label{sub-stable structure}

In this section we provide the details of the stable model category  
structure for symmetric spectra over $\C$; the result is
summarized as Theorem \ref{thm-specC}. 
We use the level model category to define the stable model category
structures on $\specC$.  
The stable model category is more difficult to
establish than the level model category, and we need to assume that $\C$ 
is a simplicial, cofibrantly generated, proper, stable model category.
The proof of the stable model structure for $\specC$
is similar to the proof of the stable model
structure for $\spec$ in~\cite[3.4]{hss},
except for one point in the proof of Proposition \ref{J-inj&W} 
where we use the stability of $\C$ instead of the fact that 
fiber sequences and cofiber sequences of spaces are stably equivalent.

Categories of symmetric spectrum objects over a model category 
have been considered more generally by Hovey in \cite{hovey-stabilization}.
Hovey relies on the general localization machinery of 
\cite{hirschhorn-book}.
Theorem \ref{thm-specC} below should be compared to 
\cite[Thms. 8.11 and 9.1]{hovey-stabilization} which are more
general but have slightly different technical assumptions.

\begin{definition}\label{stable eq.}
{\em  Let $\lambda:F_1S^1\to F_0S^0\iso S$ be the stable equivalence of
symmetric spectra which is adjoint to the identity map on the first level.  
A spectrum $Z$ in $\specC$ is an {\em $\Omega$-spectrum}
if $Z$ is fibrant on each level and the map $Z\iso Z^{F_0S^0}\to Z^{F_1S^1}$ 
induced by $\lambda$ is a level equivalence.

A map $g:A\to B$ in $\specC$ is a {\em stable equivalence} if the induced map
\[\Hom_{\specC}(g^c,Z) \ : \ \Hom_{\specC}(A^c,Z) \ \varrow{1cm} \
\Hom_{\specC}(B^c,Z) \]
is a level equivalence of symmetric spectra for any $\Omega$-spectrum $Z$;
here $(-)^c$ denotes a cofibrant replacement functor in the level
model category structure.
A map is a {\em stable cofibration} if it has the left lifting property 
with respect to each level trivial fibration, i.e., 
if it is a cofibration in the level model category structure.  
A map is a {\em stable fibration} if it has the right lifting property
with respect to each map which is both a stable cofibration and a
stable equivalence.}
\end{definition}

The above definition of $\Omega$-spectrum is just a rewrite of the usual one 
since the $n$-th level of the spectrum $Z^{F_1S^1}$ is isomorphic to 
$\Omega Z_{n+1}$. This form is more convenient here, though. 
Lemma \ref{lemma-Omega spectra} below, combined with the fact that
the stable fibrations are the $J$-injective maps shows 
that the stably fibrant objects are precisely the $\Omega$-spectra.

\begin{theorem} \label{thm-specC}
Let $\C$ be a simplicial, cofibrantly generated, proper, stable  
model category.
Then $\specC$ supports the structure of a spectral model category  
-- referred to as the {\em stable model structure} -- such that the  
adjoint functors $\Sigma^{\infty}$ and evaluation $\Ev_0$,
are a Quillen equivalence between $\C$ 
and $\specC$ with the stable model
structure.
\end{theorem}

We deduce the theorem about the stable model structure on $\specC$
from a sequences of lemmas and propositions.

\begin{lemma}\label{lemma-Omega to Omega}
Let $K$ be a cofibrant symmetric spectrum, $A$ a cofibrant spectrum
in $\specC$ and $Z$ an $\Omega$-spectrum in $\specC$.
Then the symmetric function spectrum $\Hom_{\specC}(A,Z)$  is a 
symmetric $\Omega$-spectrum and the function spectrum $Z^K$ is 
an $\Omega$-spectrum in $\specC$.
\end{lemma}
\begin{proof} By the adjunctions between smash products and
function spectra (see \ref{spectral adjunction}) we can rewrite 
the symmetric function spectrum  
$\Hom_{\specC}(A,Z)^{F_1S^1}$ as $\Hom_{\specC}(A,Z^{F_1S^1})$ in such a
way that the map $\Hom_{\specC}(A,Z)^{\lambda}$ is isomorphic to
the map $\Hom_{\specC}(A,Z^{\lambda})$; since 
\[ Z^{\lambda} \ : \  Z\iso Z^{F_0S^0}\to Z^{F_1S^1} \] 
is a level equivalence between level fibrant objects, 
the first claim follows from property 
(SP$\rscript{lv}$) of Proposition \ref{prop-level-model}.

Similarly we can rewrite the spectrum 
$(Z^K)^{F_1S^1}$ as $(Z^{F_1S^1})^K$ in such a way that the map 
$(Z^K)^{\lambda}$ is isomorphic to the map $(Z^{\lambda})^K$; 
since $Z^{\lambda}$ is a level equivalence
between level fibrant objects, the second claim follows 
from the adjoint form (SP$\rscript{lv}$(a)) 
of property (SP$\rscript{lv}$) of Proposition \ref{prop-level-model},
see~\cite[II.2 SM7(a)]{Q}.
\end{proof}

\begin{lemma} \label{lemma-stable characterization}
Let $\C$ be a simplicial, cofibrantly generated 
and left proper model category. Then a cofibration $A\to B$
is a stable equivalence if and only if 
for every $\Omega$-spectrum $Z$  the symmetric function  spectrum 
$\Hom_{\specC}(B/A,Z)$ is level contractible.
\end{lemma}
\begin{proof}
Choose a factorization of the functorial level cofibrant replacement
$A^c\to B^c$ of the given cofibration as a cofibration 
$i:A^c\to\bar B$ followed by a level equivalence $q:\bar B\to B^c$.
Then $q$ is a level equivalence between cofibrant objects, so for
every $\Omega$-spectrum $Z$, the induced map $\Hom_{\specC}(q,Z)$
is a level equivalence. 
Hence $f$ is a stable equivalence if and only if 
\[ \Hom_{\specC}(i,Z) \ : \ \Hom_{\specC}(A^c,Z) \ \varrow{1cm} \
\Hom_{\specC}(\bar B,Z) \]
is a level equivalence for every $\Omega$-spectrum $Z$. 

The symmetric spectrum $\Hom_{\specC}(\bar B/A^c,Z)$ is the fiber 
of the level fibration  $\Hom_{\specC}(i,Z)$ (by (SP$\rscript{lv}$)) 
between symmetric $\Omega$-spectra (by Lemma \ref{lemma-Omega to Omega}).
Hence the given map is a stable equivalence if and only if 
$\Hom_{\specC}(\bar B/A^c,Z)$ is level contractible.

Since $A^c\to A$ and $\bar B\to B$ are level equivalences and
$\C$ is left proper, the induced map on the cofibers
$\bar B /A^c\to B/A$ is a level equivalence between level cofibrant objects.
So for every $\Omega$-spectrum $Z$ the induced map 
$\Hom_{\specC}(B/A,Z) \to \Hom_{\specC}(\bar B/A^c,Z)$ 
is a level equivalence
(by (SP$\rscript{lv}$)) between symmetric $\Omega$-spectra 
(by Lemma \ref{lemma-Omega to Omega}).
Hence the given map is a stable equivalence if and only if 
$\Hom_{\specC}(B/A,Z)$ is level contractible, which proves the lemma.
\end{proof}

We now show that $\specC$ satisfies (SPb) of Section~\ref{sub-prerequisites}, 
an adjoint form of (SP) from Definition~\ref{def-spectral model cat}.   
This shows that $\specC$ is a spectral model category as soon 
as the stable model structure on $\specC$ is established.

\begin{proposition} \label{prop-SP}  
Let $\C$ be a simplicial, cofibrantly generated 
and left proper model category.
Let $i:A\to B$ be a cofibration in $\specC$ and 
$j:K\to L$ a stable cofibration  of symmetric spectra.
Then the pushout product map
\[ j \boxprod i \ : \ L \sm A \cup_{K \sm A}  K \sm B \ \varrow{1cm} 
\ L \sm B \]
is a cofibration in $\specC$; 
the pushout product map is a stable equivalence if  
in addition $i$ is a stable equivalence in $\specC$ or $j$ is a  
stable equivalence of symmetric spectra.
\end{proposition}

\begin{proof}
Since the cofibrations coincide in the level and the stable 
model structures for $\specC$ and for symmetric spectra, we know by property
(SP$\rscript{lv}$) of Proposition  \ref{prop-level-model} 
that $j \boxprod i$ is again a cofibration in $\specC$.  
Now suppose that one of the maps is in addition a stable equivalence.
The pushout product map $j\boxprod i$ is a cofibration with cofiber
isomorphic to $(L/K)\sm (B/A)$.
So by Lemma \ref{lemma-stable characterization} it suffices to show 
that $\Hom_{\specC}((L/K)\sm (B/A),Z)$ is level contractible 
for every $\Omega$-spectrum $Z$. 
If $i$ is a stable acyclic cofibration, then we can 
rewrite this function spectrum as
\[ \Hom_{\specC}((L/K)\sm (B/A),Z) \ \iso \ 
\Hom_{\specC}(B/A,Z^{(L/K)}) \ ; \] 
the latter spectrum is level contractible by 
Lemma \ref{lemma-stable characterization} since $Z^{(L/K)}$ 
is an $\Omega$-spectrum  by Lemma \ref{lemma-Omega to Omega}
and $i$ is a cofibration and stable equivalence. 
If $j$ is a stable acyclic cofibration, then we similarly rewrite 
the spectrum as
\[ \Hom_{\specC}((L/K)\sm (B/A),Z) \ \iso \ 
\Hom_{\spec}(L/K,\Hom_{\specC}(B/A,Z)) \ ; \]
the latter spectrum is level contractible by \cite[5.3.9]{hss}
since $\Hom_{\specC}(B/A,Z)$ 
is a symmetric $\Omega$-spectrum by Lemma \ref{lemma-Omega to Omega} 
and $L/K$ is stably contractible. 
\end{proof}

We use Theorem~\ref{thm-lifting model structure} to verify the stable 
model category structure on $\specC$.
We first define two sets $I$ and $J$ of maps in $\specC$ which will be 
generating sets for  the cofibrations and stable trivial cofibrations. 
Since the stable cofibrations are the same class of maps as the
cofibrations in the level model structure we let $I$ be the generating set 
$FI_{\C}$ which was used in Proposition \ref{prop-level-model} 
to construct the level model structure. With this choice the $I$-injectives are
precisely the level trivial fibrations.
  
The generating set for the stable trivial cofibrations is the union 
$J = FJ_{\C} \cup K$, where $FJ_{\C}$ is the generating set 
of trivial cofibrations for the level model category 
(see the proof of Proposition \ref{prop-level-model}) and $K$ is defined 
as follows. In the category of symmetric spectra over simplicial sets 
there is a map $\lambda : F_1S^1 \to  F_0S^0=S$ which is adjoint 
to the identity map on the first level; 
this map was also used in defining an $\Omega$-spectrum 
in Definition~\ref{stable eq.}.
Let $M\lambda$ be the mapping cylinder of this map, formed by taking 
the mapping cylinder of simplicial sets on each level. 
So $\lambda = r\kappa$ with $\kappa: F_1S^1 \to M\lambda$ a stable equivalence
and stable cofibration and $r: M\lambda \to S$ a simplicial 
homotopy equivalence, see \cite[3.4.9]{hss}. 
Then $K$ is the set of maps
\[ K \ = \{ \kappa\boxprod FI_{\C} \} \ = \ \{\kappa\boxprod F_ni  | i \in  
I_{\C}\} \ , \]
where for $i:A \to B$,
\[ \kappa \boxprod F_ni : (F_1S^1 \sm F_nB)  \cup_{F_1S^1 \sm F_nA}  
(M\lambda  \sm F_nA) \to M\lambda \sm F_nB \ . \]
Here we only use the pushout product, $\boxprod$, as a convenient way
of naming these maps, see also~\cite[5.3]{hss}.
Now we can verify condition (2) of Theorem~\ref{thm-lifting model structure}.
  
\begin{proposition}\label{prop-J-cof}
Let $\C$ be a simplicial, cofibrantly generated and left proper model category.
Then every relative $J$-cell complex is an $I$-cofibration and a 
stable equivalence.
\end{proposition}
\begin{proof} All maps in $J$ are cofibrations in the level model structure 
on $\specC$ of Proposition \ref{prop-level-model}, hence the
relative $J$-cell complexes are contained in the $I$-cofibrations.

We claim that for every $J$-cofibration $A\to B$ 
and every $\Omega$-spectrum $Z$, the map
\[ \Hom_{\specC}(B,Z) \ \varrow{1cm} \ \Hom_{\specC}(A,Z) \]
is a level trivial fibration of symmetric spectra. 
Hence the fiber, the symmetric spectrum $\Hom_{\specC}(B/A,Z)$,
is level contractible  and $A\to B$ is a stable equivalence
by Lemma \ref{lemma-stable characterization}.

The property of inducing a trivial fibration after applying
$\Hom_{\specC}(-,Z)$ is closed under pushout, transfinite  composition 
and retract, so by the small object argument \ref{transfinite small object} 
it suffices to check  this for the generating maps in $J=FJ_{\C}\cup K$.
The generating cofibrations in $FJ_{\C}$ are level trivial cofibrations,
so for these the claim holds by the compatibility axiom (SP$\rscript{lv}$). 
A map in the set $K$ is of the form $\kappa \boxprod F_ni$
where $\kappa:F_1S^1 \to M\lambda$ is a stable trivial cofibration of
symmetric spectra and $F_ni$ is a cofibration in $\specC$; 
hence the map $\kappa \boxprod F_ni$ is a stable trivial cofibration
between cofibrant objects by Proposition \ref{prop-SP}.
So the induced map of symmetric spectra 
$\Hom_{\specC}(\kappa \boxprod F_ni,Z)$ is a level fibration 
(by (SP$\rscript{lv}$)) between $\Omega$-spectra 
by Lemma \ref{lemma-Omega to Omega}.
In addition the fiber of the map $\Hom_{\specC}(\kappa \boxprod F_ni,Z)$
is level contractible by Lemma \ref{lemma-stable characterization},
so the map is indeed a level trivial fibration.
\end{proof}

Before turning to property (3) of Theorem~\ref{thm-lifting model
structure}, we need the following lemma.

\begin{lemma} \label{lemma-Omega spectra} 
Let $\C$ be a simplicial, cofibrantly generated model category
and $X$ a symmetric spectrum over $\C$.  Then the map $X \to *$ 
is $J$-injective if and only if $X$ is an $\Omega$-spectrum.
\end{lemma}
\begin{proof} The maps in $FJ_{\C}$ generate the trivial cofibrations 
in the level model structure of Proposition \ref{prop-level-model}, 
so $X\to *$ is $FJ_{\C}$-injective if an only if $X$ is levelwise fibrant.
Now we assume that $X$ is levelwise fibrant and show 
the map $X\to *$ is $K$-injective if and only if $X$ is an $\Omega$-spectrum. 
By adjointness $X\to *$ is $K$-injective if and only if the map 
$X^{\kappa} : X^{M\lambda}\to X^{F_1S^1}$ is a level trivial fibration. 
The projection $r:M\lambda\to S$ is a simplicial homotopy equivalence, 
so it induces a level equivalence $X\to X^{M\lambda}$. So $X\to *$
is $K$-injective if and only if the map $X^{\lambda}$ is a level equivalence,
which precisely means that $X$ is an $\Omega$-spectrum. 
\end{proof}

\begin{proposition}\label{J-inj&W}
Let $\C$ be a simplicial, cofibrantly generated, 
right proper, stable model category. 
Then a map is $J$-injective and a stable equivalence if and only if it
is a level trivial fibration.
\end{proposition}
\begin{proof} Every level equivalence is a stable equivalence 
and the level trivial fibrations are precisely the $I$-injectives.
Since the $J$-cofibrations are contained in the $I$-cofibrations,
these $I$-injectives are also $J$-injective. 
The converse is more difficult to prove.

Since $J$ in particular contains maps of the form $F_nA \to  F_nB$ 
where $n$ runs over the natural numbers and $A\to B$ runs over 
a set of generating trivial cofibrations for $\C$, $J$-injective maps 
are level fibrations. So we show that a $J$-injective stable equivalence,
$E \to B$, is a level equivalence.  Let $F$ denote the fiber and
choose a cofibrant replacement $F\rscript{c}\to F$ in the
level model category structure. Then choose a factorization in the
level model category structure
\[\begin{diagram}
\node{F\rscript{c}} \arrow{e,V} \node{E\rscript{c}}  
\arrow{e,t}{\mbox{\scriptsize lv}\, \sim} \node{E}
\end{diagram}\]
of the composite map $F\rscript{c}\to E$ as a cofibration followed by a 
level equivalence. Since $\C$ is right proper, each level of
$F \to E \to B$ is a homotopy fibration sequence in $\C$. 
Each level of
$F\rscript{c} \to E\rscript{c} \to E\rscript{c}/F\rscript{c}$
is a homotopy cofibration sequence in $\C$
(left properness is not needed here since each object is cofibrant). 
So since $\C$ is
stable, we see that $E\rscript{c}/F\rscript{c} \to B$ is a level  
equivalence.
Thus $E\rscript{c}\to E\rscript{c}/F\rscript{c}$ is a stable equivalence.
For any $\Omega$-spectrum $Z$ there is a fiber sequence of 
symmetric $\Omega$-spectra
\[ \Hom_{\C}(E\rscript{c}/F\rscript{c},Z) \ \varrow{1cm} \ 
\Hom_{\C}(E\rscript{c},Z) \ \varrow{1cm} \ 
\Hom_{\C}(F\rscript{c},Z) \]
in which the left map is a level equivalence and the right map 
is a level fibration. Hence the symmetric spectrum $\Hom_{\C}(F\rscript{c},Z)$
is level contractible which means that $F$ is stably contractible.
  
Since $F\to *$ is the pull back of the map $E\to B$, it is a $J$-injective map.
So $F$ is an $\Omega$-spectrum by Proposition \ref{lemma-Omega spectra}. 
Since $F$ is both stably contractible and an $\Omega$-spectrum, the
spectrum $\Hom_{\C}(F\rscript{c},F)$ is level contractible, 
so $F$ is level equivalent to a point. But this means that 
$E\rscript{c} \to E\rscript{c}/F\rscript{c}$, and thus also $E\to B$ 
is a level equivalence.
\end{proof}

\begin{proof}[Proof of Theorem \ref{thm-specC}]
We apply Theorem \ref{thm-lifting model structure} to show that the 
$I$-cofibrations,  stable equivalences, and $J$-injectives form 
a cofibrantly generated model category on $\specC$. Since the
$I$-cofibrations are exactly the stable cofibrations this implies 
that the $J$-injectives are the maps with the right lifting property 
with respect to the stable trivial cofibrations, i.e., 
the stable fibrations as defined before the statement of the theorem. 
The 2-out-of-3 condition, part (1) in Theorem 
\ref{thm-lifting model structure}, is clear from the definition 
of stable equivalences. Condition (2) is verified in 
Proposition \ref{prop-J-cof} and condition (3) is verified in 
Proposition \ref{J-inj&W} (since the $I$-injectives are precisely 
the levelwise trivial fibrations). So to conclude that the sets $I$
and $J$ generate a model structure with the stable equivalences as 
weak equivalence it is enough to  verify that the domains of the 
generators are small with respect to the level cofibrations. 
This has already been checked in the proof of Proposition 
\ref{prop-level-model} for the generators in $FI_{\C}$ and $FJ_{\C}$. 
So only the generators in $K$ remain. Since a pushout of objects  which are
small is small, we only need to check that $F_kA \sm M\lambda$ and 
$F_kA \sm F_1S^1$ are small with respect to the level cofibrations. 
Here $A$ is small with respect to relative $I_{\C}$-cell complexes
and hence also the cofibrations by~\cite[14.2.14]{hirschhorn-book}. 
$F_kA$ is small with respect to relative $I$-cell complexes in $\specC$ 
by adjointness and $F_1S^1$ is small with respect to all of $\spec$. 
So by various adjunctions $F_1S^1\sm F_kA$ is small  with respect to
level cofibrations. Since $M\lambda$ is the pushout of 
small objects, similar arguments show that $F_kA \sm  M\lambda$ is also 
small with respect to the level cofibrations in $\specC$.

The spectral compatibility axiom is verified in Proposition 
\ref{prop-SP} in its adjoint form (SPb). 
Thus, 
it remains to show that the adjoint functors $\Sigma^{\infty}$ and $\Ev_0$ 
are a Quillen equivalence. The suspension spectrum functor 
$\Sigma^{\infty}$ takes (trivial) cofibrations to (trivial) cofibrations 
in the level model structure. Hence $\Sigma^{\infty}$ also 
preserves (trivial) cofibrations with respect to  the stable model
structure. So the adjoint functors $\Ev_0$ and $\Sigma^{\infty}$ 
are a Quillen pair between $\C$ and $\specC$.

To show that the functors are a Quillen equivalence it suffices to show
(see \cite[Cor.\ 1.3.16]{hovey-book}) that $\Ev_0$ reflects
stable equivalences between stably fibrant objects and that 
for every cofibrant object $A$ of $\C$
the map $A\to \Ev_0 R(\Sigma^{\infty}A)$ is 
a weak equivalence where $R$ denotes any stably fibrant replacement 
in $\specC$.
So suppose that  $f:X \to Y$ is a map 
between $\Omega$-spectra with the property that 
$f_0: X_0 \to  Y_0$ is a weak equivalence in $\C$.
Since $X$ is an $\Omega$-spectrum, $X_0 \to \Omega^n X_n$ 
is a weak equivalence, and similarly for $Y$. 
Hence the map $\Omega^n f_n:\Omega^n X_n\to \Omega^nY_n$
is a weak equivalence in $\C$. Since $\C$ is stable,
the loop functor is a self-Quillen equivalence, so it reflects
weak equivalences between fibrant objects, and so $f_n:X_n\to Y_n$
is a weak equivalence in $\C$. Hence $f$ is a level, and thus a stable 
equivalence of spectra over $\C$.

Since $\C$ is stable the spectrum $\Sigf A$ (the fibrant replacement
of the suspension spectrum in the level model structure) 
is an $\Omega$-spectrum, and thus stably fibrant.  
Hence we may take $\Sigf A$ as the stably fibrant
replacement $R(\Sigma^{\infty} A)$, which proves that 
$A \to \Ev_0 R(\Sigma^{\infty}A)$ is a weak equivalence in $\C$.
\end{proof}

\subsection{The Quillen equivalence} \label{Quillen equivalence}

In this section we prove Theorem \ref{thm-main-gen}, 
i.e., we show that a suitable model category with a set 
of compact generators is Quillen equivalent to the modules over
the spectral endomorphism category of the generators.

In Theorem \ref{main theorem spectral version} 
we first formulate the result for spectral model categories;
this gives a more general result since the conditions 
about cofibrant generation and properness in $\C$ are not needed.
We then combine this with the fact that every suitable stable model category
is Quillen equivalent to a spectral model category to prove 
our main classification theorem.

\begin{definition} \label{def-EP}
{\em Let $\G$ be a set of objects in a spectral model category $\D$. 
We denote by $\EG$ the full spectral subcategory of $\D$ with objects $\G$,
i.e., $\EG(G,G')=\Hom_{\D}(G,G')$. We let
\[ \Hom(\G,-) \, : \,  \D \ \varrow{1cm} \ \mbox{mod-}\EG \]
denote the tautological functor given by $\Hom(\G,Y)(G) = \Hom_{\D}(G,Y)$.}
\end{definition}

We want to stress the reassuring fact that the stable equivalence type of
the spectral endomorphism category $\EG$ only depends on the weak equivalence 
types of the objects in the set $\G$, as long as these are all 
fibrant and cofibrant, see Corollary \ref{cor-ho invariance of endo rings}. 
This is not completely obvious since taking endomorphisms is not a functor. 

The earlier Definition \ref{def-endo ring spectrum}
of the endomorphism ring spectrum and 
endomorphism category of objects in a simplicial stable model category $\C$
is a special case of Definition \ref{def-EP} with $\D=\specC$ and $\G$ the
level fibrant replacements of the suspension spectra of the chosen
objects in $\C$. 
Again, if $\G= \{ G \}$ has a single element
then $\EG$ is determined by the single symmetric ring spectrum,
$\End_{\D}(G) = \Hom_{\D}(G,G)$.

\begin{definition} \label{def-spectral Quillen equivalence}
{\em Let $\C$ and $\D$ be spectral model categories. A {\em spectral Quillen
pair} is a Quillen adjoint functor pair $L:\C\to \D$ and $R:\D\to\C$
together with a natural isomorphism of symmetric homomorphism spectra
\[ \Hom_{\C}(A,RX) \ \iso \  \Hom_{\D}(LA,X) \]
which on the vertices of the 0-th level reduces to the adjunction isomorphism.
A spectral Quillen pair is a {\em spectral Quillen equivalence} if
the underlying Quillen functor pair is an ordinary Quillen equivalence.}
\end{definition}

In the terminology of \cite[Def.\ 4.2.18]{hovey-book} a spectral Quillen
pair would be called a `$\spec$-Quillen functor.' 

\begin{theorem} \label{main theorem spectral version} Let $\D$  
be a spectral model category and  $\G$ a set of of  
cofibrant and fibrant objects.\\
{\em (i)} The tautological functor
\[ \Hom(\G,-) \, : \, {\D} \ \varrow{1cm} \ \mbox{\em mod-}\EG \]
is the right adjoint of a spectral Quillen functor pair. 
The left adjoint is denoted $-\sm_{\EG}\G$.\\
{\em (ii)} If all objects in $\G$ are compact, then the total derived  
functors of $\Hom(\G,-)$ and $-\sm_{\EG}\G$ restrict to a triangulated  
equivalence between the homotopy category of $\EG$-modules and the localizing  
subcategory of $\Ho(\D)$ generated by $\G$.\\
{\em (iii)} If \ $\G$ is a set of compact generators for $\D$, then the  
adjoint functor pair $\Hom(\G,-)$ and $-\sm_{\EG}\G$ form a  
spectral Quillen equivalence.
\end{theorem}
\begin{proof}
(i) For an $\EG$-module $M$ the object $M\sm_{\EG}\G$ is given
by an enriched coend \cite[3.10]{kelly}. 
This means that $M\sm_{\EG}\G$ is the coequalizer of the two maps
\[
\bigvee_{G,G' \in \G} M(G') \sm \EG(G, G') \sm  G
\parallelarrows{1cm} \bigvee_{G \in \G} M(G) \sm G \ .
\]
One map in the diagram is induced by the evaluation map 
$\EG(G,G') \sm G\to G'$ and the
other is induced by the action map $M(G') \sm \EG(G,G') \to M(G)$.
The tautological functor Hom$(\G,-)$ preserves fibrations and
trivial fibrations by the compatibility axiom (SP) 
of Definition \ref{def-spectral model cat}, since all objects 
of $\G$ are cofibrant. 
So together with its left adjoint it forms a spectral Quillen pair.

\medskip

(ii) Since the functors $\Hom(\G,-)$ and $-\sm_{\EG}\G$ are a  
Quillen pair, they have adjoint total derived functors on the level of  
homotopy categories \cite[I.4]{Q}; we denote these derived functors  
by $\RHom(\G,-)$ and $-\sm^L_{\EG}\G$ respectively. 
The functor $-\sm^L_{\EG}\G$ commutes with suspension and 
preserves cofiber sequences, and the functor $\RHom(\G,-)$  commutes 
with taking loops and preserves fiber sequences \cite[I.4 Prop.\ 2]{Q}. 
In the homotopy category of a stable model  
category, the cofiber and fiber sequences coincide up to sign 
and they constitute the distinguished triangles. 
So both total derived functors preserve shifts and  
triangles, i.e., they are exact functors of triangulated categories.

For every $G\in\G$  the $\EG$-module $\Hom(\G,G)$ is isomorphic to 
the free module $F_G=\EG(-,G)$  by inspection 
and $F_G \sm_{\EG} \G$ is isomorphic to $G$ 
since they represent the same functor on $\D$. 
As a left adjoint, the functor $-\sm^L_{\EG} \G$ preserves coproducts. 
We claim that the right adjoint $\RHom(\G,-)$ also preserves coproducts. 
Since the free modules $F_G$ form a set of compact generators for the
category of $\EG$-modules (see Theorem~\ref{O-modules}), 
it suffices to show that for all $G\in\G$ 
and for every family $\{A_i\}_{i\in I}$ of objects of $\D$ the natural map
\[ \bigoplus_{i\in I} \, [F_G,\RHom(\G,A_i)]_*
^{\Ho(\mbox{\scriptsize mod-}\EG)} \  \iso \ 
[F_G, \coprod_{i\in I} \RHom(\G,A_i)]_*^{\Ho(\mbox{\scriptsize mod-}\EG)} \ 
\varrow{1cm} \]
\[ \hspace*{7cm} [F_G, 
\RHom(\G,\coprod_{i\in I} A_i)]_*^{\Ho(\mbox{\scriptsize mod-}\EG)} \] 
is an isomorphism. By the adjunctions and the identification 
$F_G \sm^L_{\EG} \G \iso G$ this map is isomorphic to the natural map
\[ \bigoplus_{i\in I} \, [G,A_i]_*^{\Ho(\D)} \  
\varrow{1cm} \ [G,\coprod_{i\in I} A_i]_*^{\Ho(\D)} \ .\] 
But this last map is an isomorphism since $G$ was assumed to be compact.
 
Both derived functors preserve shifts, triangles and coproducts; 
since they match up the free $\EG$-modules $F_G$ 
with the objects of $\G$, 
they restrict to adjoint functors between the localizing
subcategories generated by the free modules on the one side and the objects of 
$\G$ on the other side. 
We consider the full subcategories of those 
$M\in \Ho(\mbox{mod-}\EG)$ and $X\in\Ho(\D)$ respectively 
for which the unit of the adjunction
\[ \eta \ : \ M \ \varrow{1cm} \ \RHom(\G, M \sm^L_{\EG}  \G) \]
or the counit of the adjunction
\[ \nu \ : \RHom(\G,X) \sm^L_{\EG} X \ \varrow{1cm} \ X \]
are isomorphisms. Since both derived functors are exact and preserve 
coproducts, these are localizing subcategories.  
Since $F_G \sm^L_{\EG} \G \iso G$ and $\RHom(\G,G)\iso F_G$, 
the map $\eta$ is an isomorphism for every free module, and the map $\nu$ 
is an isomorphism for every object of $\G$.
Since the free modules $F_G$ generate the homotopy category 
of $\EG$-modules, the  claim follows.

\medskip

(iii) Now the localizing subcategory generated by $\G$ is the 
entire homotopy category of $\D$, so part (ii) of the theorem implies 
that the total derived functors of $\Hom(\G,-)$ and \mbox{$-\sm_{\EG}\G$} 
are inverse equivalences of homotopy categories. 
Hence this pair is a Quillen equivalence.
\end{proof}

Now we can finally give the

\begin{proof}[Proof of Theorem \ref{thm-main-gen}:]
We can combine Theorem \ref{thm-specC} and Theorem \ref{main theorem  
spectral version} (iii) to obtain a diagram of  model categories and  
Quillen equivalences
\vspace{-0cm}
\[
\raisebox{0.18cm}{$\C$}
\parbox{2cm}{\begin{center} \raisebox{-0.25cm}{$\scriptstyle  
\Sigma^{\infty}$} \\ \hspace*{0.7cm}
$\begin{diagram} \arrow{e} \end{diagram}$ \\ \vspace{-0.30cm}  
\hspace*{2.2cm} $\begin{diagram} \arrow{w} \end{diagram}$ \\  
\raisebox{0.4cm}{\mbox{\scriptsize Ev}$_0$} \end{center}}
\raisebox{0.18cm}{$\specC$}
\parbox{2cm}{\begin{center} \raisebox{-0.25cm}{$\scriptstyle -  
\sm_{\EG}\G$} \\  \hspace*{2.2cm}
$\begin{diagram} \arrow{w} \end{diagram}$ \\ \vspace{-0.30cm}  
\hspace*{0.7cm} $\begin{diagram} \arrow{e} \end{diagram}$ \\  
\raisebox{0.4cm}{$\scriptstyle \Hom(\G,-)$} \end{center}}
\raisebox{0.18cm}{mod-$\EG$}
\vspace{-0.5cm}
\]
(the left adjoints are on top). 
First, we may assume that each object in the set $\Pc$ of compact generators
for $\C$ is cofibrant.
Since the left Quillen functor pair above induces an equivalence 
of homotopy categories, the suspension spectra of the objects in $\Pc$ 
form a set of compact cofibrant generators for  $\specC$. 

We denote by $\G$ the set of level fibrant replacements $\Sigf P$ of 
the given generators in $\C$.
The spectral endomorphism category $\EG$ in the sense 
of Definition \ref{def-EP} is equal to the endomorphism 
category $\E(\Pc)$ associated to $\Pc$ 
by Definition \ref{def-endo ring spectrum}. 
Since $\C$ is stable, the spectra $\Sigf P$ are $\Omega$-spectra,
hence they are both fibrant (by Lemma \ref{lemma-Omega spectra}) 
and cofibrant in the stable model structure on $\specC$.
So we can apply Theorem \ref{main theorem spectral version} (iii)  
to get the second Quillen equivalence.
\end{proof}

\begin{remark}{\em {\bf Finite localization and $\EP$-modules.}
\label{rem-fin-loc} 
Suppose $\Pc$ is a set of compact objects of a triangulated category $\T$ 
with infinite coproducts. Then there always exists an idempotent 
localization functor $L_{\Pc}$ on $\T$ whose acyclics are precisely 
the objects of the localizing subcategory generated by $\Pc$ 
(compare \cite{miller-finite} or the proofs of Lemma 2.2.1 or 
\cite[Prop.\ 2.3.17]{hps}). These localizations are often referred to 
as {\em finite Bousfield localizations} away from $\Pc$. 

Theorem \ref{main theorem spectral version} gives a lift 
of finite localization to the model category level. 
Suppose $\C$ is a stable model category with a set $\Pc$ 
of compact objects, and let $L_{\Pc}$ denote the associated
localization functor on the homotopy category of $\C$. 
By Theorem 5.3.2 (ii) the acyclics for $L_{\Pc}$
are equivalent to the homotopy category of $\EP$-modules, 
the equivalence arising from a Quillen adjoint functor pair. 
Furthermore the counit of the derived adjunction
\[ \Hom(\Pc,X) \, \sm_{\EP}^L \, \Pc \ \varrow{1cm} \ X \]
is the acyclicization map and its cofiber is a model for the 
localization $L_{\Pc} X$.}
\end{remark}

\section{Morita context} \label{sec-Morita}
\stepcounter{subsection}

In the classical algebraic context there is a characterization of 
equivalences of module categories in terms of bimodules,
see for example \cite[\S 22]{anderson-fuller}.
We provide an analogous result for module categories over ring spectra. 
As usual, here instead of actual equivalences of module categories 
one obtains Quillen equivalences of model categories. 
We state the Morita context for symmetric ring spectra
and spectral Quillen equivalences
(see Definition \ref{def-spectral Quillen equivalence}).
  
\begin{definition} 
{\em If $R$ is a symmetric ring spectrum and $\C$ a spectral model category, 
then an $R$-$\C$-bimodule is an object $X$ of $\C$ on which $R$ acts 
through $\C$-morphisms, i.e., a homomorphism of
symmetric ring spectra from $R$ to the endomorphism ring spectrum of $X$.}
\end{definition}

If $\C$ is the category of modules over another symmetric ring spectrum $T$, 
then this notion of bimodule specializes to the usual one, i.e., an 
$R$-$(T\mbox{-mod})$-bimodule is the same 
as a right $R\rtiny{op}\sm T$-module. 
In the following Morita theorem, the implication $(3) \Longrightarrow (2)$ 
is a special case of the classification of monogenic stable model categories 
(Theorem \ref{thm-main-one}); hence the implication $(2) \Longrightarrow (3)$
is a partial converse to that classification result.
The implication $(2) \Longrightarrow (1)$ says that certain chains of 
Quillen equivalences can be rectified into a single 
Quillen equivalence whose left adjoint is given by smashing with a bimodule.

\begin{theorem} \label{thm-Morita} {\bf (Morita context)} 
Consider the following statements for a symmetric ring
spectrum $R$ and a spectral model category $\C$.
\renewcommand{\labelenumi}{(\arabic{enumi})}
\begin{enumerate}
\item There exists an $R$-$\C$-bimodule $M$ such that smashing with $M$ 
over $R$ is the left adjoint of a Quillen equivalence 
between the category of $R$-modules and $\C$.
\item There exists a chain of spectral Quillen equivalences 
through spectral model categories between the category of $R$-modules and $\C$.
\item The category $\C$ has a compact, cofibrant and fibrant generator $M$ 
such that $R$ is stably equivalent to the endomorphism ring spectrum of $M$.
\end{enumerate}
Then conditions {\em (2)} and {\em (3)} are equivalent, and condition {\em (1)}
implies both conditions {\em (2)} and {\em (3)}.
If $R$ is cofibrant as a symmetric ring spectrum, then all three conditions 
are equivalent.
  
Furthermore, if $\C$ is the category of modules over another 
symmetric ring spectrum $T$ which is cofibrant as a symmetric spectrum, 
then condition {\em (1)} is equivalent to the condition\\
\hspace*{0.5cm} {\em (4)} There exists an $R$-$T$-bimodule $M$ and a 
$T$-$R$-bimodule $N$ which are cofibrant as right modules such that
\begin{itemize}
\item  $M \sm_T N$ is stably equivalent to $R$ as an $R$-bimodule and
\item  $N \sm_R M$ is stably equivalent to $T$ as a $T$-bimodule.
\end{itemize}
\end{theorem}

\begin{remark} {\em The cofibrancy conditions in the Morita theorem 
can always be arranged since every symmetric ring spectrum has 
a stably equivalent cofibrant replacement in the model category 
of symmetric ring spectra \cite[5.4.3]{hss}; 
furthermore the underlying symmetric spectrum of a cofibrant ring spectrum 
is always cofibrant (\cite[4.1]{ss} or \cite[5.4.3]{hss}).
It should not be surprising that cofibrancy conditions 
have to be  imposed in the Morita theorem. 
In the algebraic context the analogous conditions show up in Rickard's paper 
\cite{rickard2}: when trying to realize derived equivalences between 
$k$-algebras by derived tensor product with bimodule complexes, 
he has to assume that the algebras are flat over the ground ring $k$, 
see \cite[Sec.\ 3]{rickard2}.}
\end{remark}

\begin{proof}[Proof of the Morita theorem]
Condition (1) is a special case of (2).
 
$(2) \Longrightarrow (3)$: 
This implication follows from the homotopy invariance of endomorphism
ring spectra under spectral Quillen equivalences, 
see  Corollary \ref{cor-ho invariance of endo rings}.
We choose a chain of spectral Quillen equivalences through spectral model
categories.  Then we choose a trivial cofibration $R \to R\rscript{f}$
of $R$-modules such that $R\rscript{f}$ is fibrant;
then $R$ is stably equivalent to the endomorphism ring spectrum 
of $R\rscript{f}$. We define an object $M$ of $\C$
by iteratively applying the functors in the chain of Quillen equivalences, 
starting with $R\rscript{f}$. 
In addition we take a fibrant or cofibrant replacement after each 
application depending on whether we use a left or right Quillen functor. 
By a repeated application of Corollary \ref{cor-ho invariance of endo rings}
the endomorphism ring spectra of these objects, including the final one, $M$,
are all stably equivalent to $R$.  
By construction the object $M$ is isomorphic in the homotopy category of $\C$ 
to the image of the free $R$-module of rank one under the equivalence 
of homotopy categories induced by the Quillen equivalences.
Hence $M$ is also a compact generator for $\C$ and satisfies condition (3).

$(3) \Longrightarrow (2)$: 
This is essentially a special case of 
Theorem \ref{main theorem spectral version} (iii).  More  precisely, that
theorem constructs a spectral Quillen equivalence between $\C$ 
and the category of modules over the endomorphism ring spectrum 
of the generator $M$. Furthermore, restriction and extension
of scalars are spectral Quillen equivalences for two stably equivalent 
ring spectra \cite[Thm.\ 5.4.5]{hss}, see also~\cite[Thm.\ 11.1]{mmss}
or Theorem~\ref{O-modules}, which establishes condition (2).

$(3) \Longrightarrow (1)$, provided $R$ is cofibrant 
as a symmetric ring spectrum: 
Since $M$ is  fibrant and cofibrant the endomorphism ring spectrum 
$\End_{\C}(M)$ is fibrant. Since $R$ is cofibrant as a symmetric
ring spectrum any isomorphism between $R$ and $\End_{\C}(M)$ 
in the homotopy category of symmetric ring spectra can be lifted 
to a stable equivalence $\eta : R \to \End_{\C}(M)$. 
In particular the stable equivalence $\eta$ makes $M$ into an 
$R$-$\C$-bimodule. The functor $X \mapsto X \sm_R M$ is left adjoint to the
functor $Y \mapsto \Hom_{\C}(M,Y)$ from $\C$ to the category of $R$-modules. 
To show that these adjoint functors form a Quillen equivalence we note 
that they factor as the composite of two adjoint functor pairs
\vspace{-0cm}
\[
\raisebox{0.18cm}{$\C$}
\parbox{2cm}{\begin{center} 
\raisebox{-0.25cm}{$\scriptstyle  -\sm_{\End_{\C}(M)}M$} \\ 
\hspace*{2.2cm}
$\begin{diagram} \arrow{w} \end{diagram}$ \\ \vspace{-0.30cm}  
\hspace*{0.7cm} $\begin{diagram} \arrow{e} \end{diagram}$ \\  
\raisebox{0.4cm}{$\scriptstyle \Hom_{\C}(M,-)$} \end{center}}
\raisebox{0.18cm}{mod-$\End_{\C}(M)$}
\parbox{2cm}{\begin{center} \raisebox{-0.25cm}{$\scriptstyle \eta_*$} \\  \hspace*{2.2cm}
$\begin{diagram} \arrow{w} \end{diagram}$ \\ \vspace{-0.30cm}  
\hspace*{0.7cm} $\begin{diagram} \arrow{e} \end{diagram}$ \\  
\raisebox{0.4cm}{$\scriptstyle \eta^*$} \end{center}}
\raisebox{0.18cm}{mod-$R$}
\vspace{-0.5cm}
\]
(the left adjoints are on top). Since $M$ is a cofibrant and fibrant 
compact generator for $\C$, the left pair is a Quillen equivalence by 
Theorem  \ref{main theorem spectral version} (iii). The other adjoint functor
pair is restriction and extension of scalars along the stable equivalence 
of ring spectra $\eta : R \to  \End_{\C}(M)$, which is a Quillen equivalence 
by \cite[Thm.\ 5.4.5]{hss} or~\cite[Thm.\ 11.1]{mmss}.

For the rest of the proof we assume that $\C$ is the category of modules over 
a symmetric ring spectrum $T$ which is cofibrant as a symmetric spectrum.\\
$(1) \Longrightarrow (4)$: 
Since smashing with $M$ over $R$ is a left Quillen functor and
since the free $R$-module of rank one is cofibrant, $M \iso R \sm_RM$ 
is cofibrant as a right $T$-module.
We choose a fibrant replacement $T\to T\rscript{f}$ of $T$ in the category 
of $T$-bimodules and we let $N$ be a cofibrant replacement 
of the $T$-$R$-bimodule $\Hom_T(M,T\rscript{f})$. The forgetful functor from
$T$-$R$-bimodules to right $R$-modules preserves cofibrations since 
its right adjoint
\[ \Hom_{\spec}(T,-) : \mbox{mod-}R \ \varrow{1cm} \  T\mbox{-mod-}R \]
preserves trivial fibrations (because $T$ is cofibrant 
as a symmetric spectrum).  In particular the
bimodule $N$ is cofibrant as a right $R$-module. We will exhibit two chains 
of stable equivalences of bimodules
\[  N \sm_R M \ \varr{1cm}{\sim }\  T\rscript{f} \ \varl{1cm}{\sim} \ T  
\mbox{\quad  and \quad} M \sm_T N \ \varr{1cm}{\sim} \ 
\Hom_T(M,(M \sm_T T\rscript{f})\rscript{f}) \ 
\varl{1cm}{\sim} \  R \]
where $M \sm_T T\rscript{f} \to  (M \sm_T T\rscript{f})\rscript{f}$ 
is a fibrant approximation in the category of $R$-$T$-bimodules. 
This will establish condition (4).
  
Since $-\sm_R M$ was assumed to be a left Quillen equivalence 
and the approximation map $N \to \Hom_T(M,T\rscript{f})$ 
is a weak equivalence, 
so is its adjoint $N \sm_RM \to T\rscript{f}$; but this  adjoint is
even a map of $T$-bimodules. The equivalence $T\to T\rscript{f}$ was chosen in 
the beginning. For the next equivalence we smash the $T$-bimodule equivalence 
$N \sm_R M \to T\rscript{f}$ with the right-cofibrant
bimodule $M$ to get a stable equivalence of $R$-$T$-bimodules
\[ M \sm_TN \sm_R M \ \varr{1cm}{\sim} \ M \sm_T T\rscript{f}  \ . \]
We then compose with the approximation map and obtain a stable equivalence 
of $R$-$T$-bimodules $M \sm_T N \sm_R M \to (M \sm_T T\rscript{f})\rscript{f}$.
Since $M$ and $N$ are cofibrant as right modules, the $R$-bimodule
$M \sm_T N$ is cofibrant as a right $R$-module. Since $-\sm_RM$ 
is a left Quillen equivalence, the adjoint
$M \sm_T N \to \Hom_T(M,(M \sm_T T\rscript{f})\rscript{f})$ 
is thus a stable equivalence of $R$-bimodules. 
For the same reason the adjoint of the composite stable equivalence
\[  R \sm_RM  \ \iso \ M\sm_T T \ \varr{1cm}{\sim} \ M \sm_T T\rscript{f} \ 
\varr{1cm}{\sim} \ ( M \sm_T T\rscript{f})\rscript{f} \]
gives the final stable equivalence of $R$-bimodules 
$R \to \Hom_T(M,(M \sm_T T\rscript{f})\rscript{f})$.

$(4) \Longrightarrow (1)$: Let $M$ and $N$ be bimodules which satisfy 
the conditions of (4). The functor $X \mapsto X \sm_RM$ is left adjoint 
to the functor $Y \mapsto \Hom_T(M,Y)$ from the category of $T$-modules to
the category of $R$-modules. Since $M$ is cofibrant as a right $T$-module, 
this right adjoint preserves fibrations and trivial fibrations 
by the spectral axiom (SP). So $-\sm_R M$ and $\Hom_T(M,-)$
form a Quillen functor pair. By condition (4) the left derived functor 
of $-\sm_R M$  is an equivalence of derived categories 
(with inverse the left derived functor of smashing with the bimodule $N$).
So the functor $ -\sm_R M$ is indeed the left adjoint of a Quillen equivalence.
\end{proof}

\section{A generalized tilting theorem}
\label{sec-tilting} \stepcounter{subsection}

In this section we state and prove a generalization of Rickard's 
``Morita theory for derived categories'', \cite{rickard1}. 
Rickard studies the question of when two rings  
are derived equivalent, i.e., when there exists a triangulated equivalence 
between various  derived categories of the module categories. 
He shows \cite[Thm.\ 6.4]{rickard1} that a necessary and sufficient condition
for such a derived equivalence is the existence of a tilting complex. 
A tilting complex for a pair of rings $\Gamma$ and $\Lambda$ is 
a bounded complex $X$ 
of finitely generated projective $\Gamma$-modules which generates
the derived category and whose graded ring of self extension groups 
$[X,X]^{\D(\Gamma)}_*$ is isomorphic to $\Lambda$, concentrated in dimension 
zero.

 We generalize the result in two directions. 
First, we allow the input to be a stable model category (which generalizes
categories of chain complexes of modules). Second, we allow
for a set of special generators, as opposed to a single tilting complex.
The compact objects in the unbounded derived category of a ring are precisely 
the perfect complexes, i.e., those chain complexes which are 
quasi-isomorphic to a bounded complex of finitely
generated projective modules \cite[Prop.\ 6.4]{boek-nee}.
In our context we thus define a set of {\em tiltors} in
a stable model category to be a set $\mathbb T$ of compact generators such that
for any two objects $T,T'\in \mathbb T$ the graded homomorphism group
$[T,T']_*^{\HoC}$ in the homotopy category is concentrated in dimension zero. 
Then Theorem \ref{thm-tilting} shows that the existence of a set of
tiltors is necessary and sufficient for a stable model category 
to be Quillen equivalent or derived equivalent to the category of
chain complexes over a {\em ringoid} (ring with several objects).
Recall that a ringoid is a small category whose hom-sets carry an
abelian group structure for which composition is bilinear. 
A {\em module} over a ringoid is a contravariant
additive functor to the category of abelian groups.

Of course an interesting special case is that of a stable
model category with a single tiltor, i.e., a single compact generator
whose graded endomorphism ring in the homotopy category is concentrated
in dimension zero. Then ringoids simplify to rings.
In particular when $\C$ is the model category of chain complexes 
of $\Gamma$-modules for some ring $\Gamma$, then a single tiltor is 
the same (up to quasi-isomorphism) as a tilting complex, and
the equivalence of conditions (2') and (3) below recovers Rickard`s 
`Morita theory for derived categories' \cite[Thm.\ 6.4]{rickard1}.  

Condition (1) in the tilting theorem refers to a standard model structure
on the category of chain complexes of  $\uA$-modules.
The model structure we have in mind is the {\em projective} model
structure: the weak equivalences are the quasi-isomorphisms 
and the fibrations are the epimorphisms.
Every cofibrations in this model structure is in particular a monomorphism
with degreewise projective cokernel, but for unbounded complexes 
this condition is not sufficient to characterize a cofibration.
In the single object case, i.e., for modules over a ring, the projective
model structure on complexes is established in~\cite[Thm.\ 2.3.11]{hovey-book}.
For modules over a ringoid the arguments are very similar, just that the
free module of rank one has to be replaced by the set of free
(or representable) modules $F_a=\uA(-,a)$ for $a\in\uA$.
This model structure can be established using 
a version of Theorem~\ref{O-modules} where the enrichment over 
$\spec$ is replaced by one over chain complexes.
The projective model structure for complexes of $\uA$-modules is also a
special case of \cite[Thm.~5.1]{Christensen-Hovey-relative}.
Indeed, the projective (in the usual sense) $\uA$-modules together
with the epimorphisms form a projective class 
(in the sense of \cite[Def.~1.1]{Christensen-Hovey-relative}),
and this class is determined
(in the sense of \cite[Sec.~5.2]{Christensen-Hovey-relative})
by the set of small, free modules $\{F_a\}_{a\in\uA}$.

\begin{theorem} \label{thm-tilting} {\bf (Tilting theorem)}
Let $\C$ be a simplicial, cofibrantly generated, proper, stable
model category and  $\uA$ a ringoid.
 Then the following conditions are equivalent:\\
\hspace*{0.5cm}{\em (1)} There is a chain of Quillen equivalences between $\C$ 
and the model category of chain complexes of  $\uA$-modules.\\
\hspace*{0.5cm}{\em (2)} The homotopy category of $\C$ is 
triangulated equivalent to $\D(\uA)$, the unbounded derived
category of the ringoid $\uA$.\\
\hspace*{0.5cm}{\em (2')}  $\C$ has a set of compact generators 
and the full subcategory of compact objects in $\HoC$ is
triangulated equivalent to $K^b(\mbox{\em proj-}\uA)$, 
the homotopy category of bounded chain complexes of
finitely generated projective  $\uA$-modules.\\
\hspace*{0.5cm}{\em (3)} The model category $\C$ has a set of tiltors
whose endomorphism ringoid in the  homotopy category of $\C$ 
is isomorphic to $\uA$.
\end{theorem}

\begin{example} \label{ex-rat equivariant via tilting}
{\em Let $G$ be a finite group. As in 
Example \ref{multiple generator examples} (i) the category
of $G$-equivariant orthogonal spectra~\cite{mm} 
based on a complete universe $\U$ is a simplicial, stable 
model category, and the equivariant suspension spectra of the homogeneous 
spaces $G/H^+$ form a set of compact generators as $H$ runs through 
the subgroups of $G$. 
Rationalization is a smashing Bousfield localization so 
the rationalized suspension spectra form a set of compact generators 
of the rational $G$-equivariant stable homotopy category. 
The homotopy groups of the
function spectra between the various generators are torsion in dimensions
different from zero \cite[Prop.\ A.3]{greenlees-may}, 
so the rationalized suspension spectra form a set of tiltors.
Modules over the associated ringoid are nothing but rational Mackey
functors, so the tilting theorem \ref{thm-tilting} shows that the rational 
$G$-equivariant stable homotopy category
is equivalent to the derived category of rational Mackey functors.
In turn, since these rational Mackey functors are all projective and injective,
the derived category is equivalent to the graded category.  So this recovers
the Theorem A.1 in \cite{greenlees-may}. 

For a non-complete universe, one considers rational $U$-Mackey 
functors~\cite{gaunce} which are modules over the endomorphism ringoid of the
rationalized suspension spectra of $G/H^+$. For example, for the trivial 
universe $U$, these rational $U$-Mackey functors are rational coefficient 
systems. The rational $G$-equivariant stable homotopy category based 
on a non-complete universe $U$ is then equivalent to the derived category 
of the associated rational $U$-Mackey functors.
} 
\end{example}

\begin{example} \label{ex-pure via tilting}
{\em Let $A$ be a ring and consider the {\em pure projective} model category 
structure in the sense of Christensen and 
Hovey~\cite[5.3]{Christensen-Hovey-relative}
on the category of chain complexes of $A$-modules 
(see also Example \ref{stable examples} (xiii)).
A map $X\to Y$ of complexes is a weak equivalence if and only if for
every finitely generated $A$-module $M$ the induced map of mapping complexes
$\Hom_A(M,X)\to\Hom_A(M,Y)$ is a quasi-isomorphism. 
Let $\G$ be a set of representatives of the isomorphism classes of
indecomposable finitely generated $A$-modules. 
Then $\G$ forms a set of compact generators for the pure derived category 
$\D_{\Pc}(A)$. Since furthermore every finitely generated module is pure
projective, maps in the pure derived category between modules in $\G$ 
are concentrated in dimension zero. In other words, the 
indecomposable finitely generated $A$-modules form a set of tiltors.
So Theorem \ref{thm-tilting} implies that the pure projective model category
of $A$ is Quillen equivalent to the modules over the ringoid given
by the full subcategory of $A$-modules with objects $\G$.}
\end{example}

\begin{remark} {\em We want to emphasize one special feature 
of the tilting situation.  For general
stable model categories the notion of Quillen equivalence 
is considerably stronger than triangulated equivalence of homotopy categories 
(see Remark \ref{rem-K(n)} for an example).  Hence it is somewhat remarkable 
that for chain complexes of modules over ringoids the two notions
are in fact equivalent. In general the homotopy category determines 
the homotopy groups of the spectral endomorphism category, 
but not its homotopy type. The real  reason behind
the equivalences of conditions (1) and (2) above is the fact that 
in contrast to arbitrary ring spectra or spectral categories, 
Eilenberg-Mac Lane objects are determined by their homotopy groups,
see Proposition \ref{prop-EMuniqueness}.}
\end{remark}

As a tool for proving the generalized tilting theorem we introduce the
{\em Eilenberg-Mac Lane spectral category} $H\uA$ of a ringoid $\uA$.
This is simply the many-generator version of the symmetric Eilenberg-Mac Lane
ring spectrum~\cite[1.2.5]{hss}. 
The key property is that module spectra over the 
Eilenberg-Mac Lane spectral category $H\uA$ are Quillen equivalent to
chain complexes of $\uA$-modules. 

\begin{definition} \label{def-EM spectral category} {\em
Let $\uA$ be a ringoid. The Eilenberg-Mac Lane spectral category $H\uA$
is defined by 
\[ H\uA \ = \ \uA \tensor H\mZ \ , \] 
where $H\mZ$ is the symmetric Eilenberg-Mac Lane 
ring spectrum of the integers~\cite[1.2.5]{hss}. 
In more detail, $H\uA$ has the same set of objects as $\uA$, 
and the morphism spectra are defined by
\[ H\uA(a,b)_p \ = \ \uA(a,b)\tensor \widetilde{\mathbb Z}[S^p] \ . \]
Here $\widetilde{\mathbb Z}[S^p]$ denotes the reduced simplicial
free abelian group generated by the pointed simplicial set 
$S^p= S^1\sm\dots\sm S^1$ ($p$ factors), and the symmetric group permutes the
factors. Composition is given by the composite 
\begin{eqnarray*}
 H\uA(b,c)_p \sm  H\uA(a,b)_q & = & 
(\uA(b,c)\tensor \widetilde{\mathbb Z}[S^p]) \sm  
(\uA(a,b)\tensor \widetilde{\mathbb Z}[S^q]) \\
& \varr{1.5cm}{\text{shuffle}} &
\uA(b,c)\tensor\uA(a,b) \tensor 
\widetilde{\mathbb Z}[S^p]\tensor \widetilde{\mathbb Z}[S^q] \\
& \varr{1.5cm}{\circ \, \tensor\, \iso} &
\uA(a,c) \tensor \widetilde{\mathbb Z}[S^{p+q}] \ = \ H\uA(a,c)_{p+q}  \ . 
\end{eqnarray*}
The unit map
\[  S^p \ \varrow{1cm} \ \uA(a,a)\tensor \widetilde{\mathbb Z}[S^p] \ = \
H\uA(a,a)_p \] 
is the inclusion of generators.
}
\end{definition}

We prove the following result in Appendix \ref{app-EM}.

\begin{theorem}\label{thm-chains and EM} For any ringoid $\uA$,  
the category of complexes of $\uA$-modules and the category 
of modules over the Eilenberg-Mac Lane spectral category $H\uA$
are Quillen equivalent,
\[ \mbox{\em mod-}H\uA \simeq_Q \mbox{\em Ch}_{\uA} \]
\end{theorem}

\begin{proof}[Proof of Theorem \ref{thm-tilting}]
Every Quillen equivalence of stable model categories induces an equivalence of 
triangulated homotopy categories, so condition (1) implies condition (2). 
Any triangulated equivalence restricts to an equivalence 
between the respective subcategories of compact objects. 
By the same argument as \cite[Prop.\ 6.4]{boek-nee} (which deals with
the special case of complexes of modules over a ring),
$K^b(\mbox{proj-}\uA)$ is equivalent to the full subcategory
of compact objects in $\D(\uA)$.
Since the derived category of a ringoid has a set of compact generators, 
so does any equivalent triangulated category. 
Hence condition (2) implies condition (2').
  
Now we assume condition (2') and we choose a triangulated equivalence 
between $K^b(\mbox{proj-}\uA)$ and the full subcategory of compact
objects in $\HoC$. 
For $a\in \uA$ we let $T_a$ be a representative in $\C$ of the image 
of the representable $\uA$-module $F_a=\uA(-,a)$, viewed as a complex 
concentrated in dimension zero.
Since the collection of modules $\{F_a\}_{a\in\uA}$ is a set of tiltors
for the derived category, the set $\{T_a\}_{a\in\uA}$ has
all the properties of a set of tiltors, 
except that it may not generate the full homotopy category of $\C$. 
However the localizing subcategory generated by the  $T_a$'s coincides 
with the localizing subcategory generated by all compact objects 
since on the other side of the equivalence the complexes $F_a$ 
generate the category $K^b(\mbox{proj-}\uA)$. In general the compact objects 
might not generate all of $\HoC$ (see \cite[Cor.\ B.13]{hovey-strickland-kn} 
for some extreme cases where the zero object is the only compact object),
but here this is assumed in (2').
So the $T_a$'s generate $\C$, hence they are a set of tiltors, 
and so condition (3) holds.
  
If on the other hand $\C$ has a set of tiltors $\mathbb T$,
then $\mathbb T$ is in particular a set of compact generators, 
so by Theorem \ref{thm-main-one}, $\C$ is 
Quillen equivalent to the category of modules over the endomorphism
category $\End(\mathbb T)$. 
In this special case the homotopy type of the spectral category 
$\End(\mathbb T)$ is determined by its homotopy groups: 
since the homotopy groups of $\End(\mathbb T)$ are concentrated in dimension 0,
$\End(\mathbb T)$ is stably equivalent to $H\uA$, the 
Eilenberg-Mac Lane spectral category of its component ringoid $\uA$
by Proposition \ref{prop-EMuniqueness}. Thus the categories of 
$\End(\mathbb T)$-modules and $H\uA$-modules are Quillen equivalent 
by Theorem \ref{O-modules}. 
Theorem \ref{thm-chains and EM} gives the final step 
in the chain of Quillen-equivalences
\[ \C \quad \simeq_Q \quad \text{mod-}\End(\mathbb T) \quad \simeq_Q \quad 
\text{mod-}H\uA \quad \simeq_Q \quad \mbox{Ch}\uA \ . \] 
\end{proof}

\begin{appendix}
\section{Spectral categories}

In this appendix we develop some general theory of
modules over spectral categories (Definition \ref{def-spectral cat}). 
The arguments are very similar to the case of spectral categories 
with one object, i.e., symmetric ring spectra. 

\subsection{Model structures for modules over spectral categories}
\label{subsec-O-mod}

A {\em morphism} $\Psi:\Oc\to \R$ of spectral categories is simply
a spectral functor. The {\em restriction of scalars}
\[ \Psi^* \ : \ \mbox{mod-}\R \ \varrow{1cm} \ \mbox{mod-}\Oc \quad , 
\quad M \ \longmapsto \ M \circ \Psi \]
has a left adjoint functor $\Psi_*$, also denoted $-\sm_{\Oc} \R$,
which we refer to as {\em extension of scalars}. 
As usual it is given by an enriched coend,
i.e., for an \Omod\ $N$ the $\R$-module $\Psi_* N = N\sm_{\Oc} \R$ 
is given by the coequalizer of the two $\R$-module homomorphisms 
\[ \bigvee_{o,p\in \Oc} N(p) \, \sm \, \Oc(o,p)\, \sm \, F_{\Psi(o)} 
\parallelarrows{1cm} \bigvee_{o\in\Oc} N(o)\, \sm\, F_{\Psi(o)}  \ , \] 
where $F_{\Psi(o)}=\R(-,\Psi(o))$ is the free $\R$-module associated
to the object $\Psi(o)$.
We call $\Psi:\Oc\to\R$ a {\em stable equivalence} of spectral categories
if it is a bijection on objects and if for all objects $o,o'$ in $\Oc$
the map
\[ \Psi_{o,o'} \ : \ \Oc(o,o') \ \varrow{1cm} \ \R(\Psi(o),\Psi(o')) \]
is a stable equivalence of symmetric spectra.

Next we establish the model category structure for $\Oc$-modules, 
show its invariance
under restriction of scalars along a stable equivalence of spectral categories
and exhibit a set of compact generators.

\begin{theorem} \label{O-modules}
\begin{enumerate}
\item Let $\Oc$ be a  spectral category.
Then the category of \Omods\ with the objectwise stable equivalences,
objectwise stable fibrations, and cofibrations is a cofibrantly generated
spectral model category. 
\item The free modules $\{F_o\}_{o\in\Oc}$ given by
$F_o=\Oc(-,o)$ form a set of compact generators 
for the homotopy category of $\Oc$-modules.
\item Assume $\Psi:\Oc \to \R$ is a stable equivalence of spectral categories. 
Then restriction and extension of scalars along $\Psi$ form a spectral Quillen
equivalence of the module categories.
\end{enumerate}
\end{theorem}

\begin{proof}
We use~\cite[2.3]{ss} to lift the stable model structure from (families of)
symmetric spectra to $\Oc$-modules. 
Let $S_{\Oc}$ denote the spectral category with the same set of objects
as $\Oc$, but with morphism spectra given by
\[ S_{\Oc}(o,o') \ = \ \left\lbrace \begin{array}{ll}
S & \mbox{if $o=o'$,} \\ * & \mbox{else.}
\end{array} \right. \]
An $S_{\Oc}$-module is just a family of symmetric spectra 
indexed by the objects of $\Oc$.  
Hence it has a cofibrantly generated model category structure
in which the cofibrations, fibrations and weak equivalences are defined
objectwise on underlying symmetric spectra.  
Here the generating trivial cofibrations 
are maps between modules concentrated at one object, i.e. of the form $A_o$
with $A_o(o)= A$ and $A_o(o')=\ast$ if $o \ne o'$.  

The unit maps give a morphism of spectral categories $S_{\Oc} \to \Oc$,
which in turn induces adjoint functors of restriction 
and extension of scalars between the module categories.  
This produces a triple $- \sm_{S_{\Oc}} \Oc$ on $S_{\Oc}$-modules 
with the algebras over this triple the $\Oc$-modules.   
Then the generating trivial cofibrations for $\Oc$-modules are maps between 
modules of the form $A_o \sm_{S_{\Oc}} \Oc = A \sm \Oc(-, o)= A \sm F_o$.
Hence the monoid axiom~\cite[5.4.1]{hss}
applies to show that the new generating trivial cofibrations and their
relative cell morphisms are weak equivalences.  
Thus, since all symmetric spectra, hence all $S_{\Oc}$-modules 
are small, the model category structure follows by criterion (1) 
of~\cite[2.3]{ss}. We omit the verification of the spectral axiom (SP),
which then implies stability by Lemma \ref{lemma-spectral implies stable}. 

The proof of (ii) uses the adjunction defined above between $S_{\Oc}$-modules
and $\Oc$-modules.  Since $F_o= S_o \sm_{S_{\Oc}} \Oc$,
\[ [F_o, M]_*^{\Ho(\mbox{\scriptsize mod-}\Oc)} \iso 
[S_o, M]_*^{\Ho(\mbox{\scriptsize mod-}S_{\Oc})} \iso [S,M(o)]_*^{\Ho(\spec)}. \]
Thus, since $S$ is a generator for $\spec$ and an $\Oc$-module is trivial
if and only if it is objectwise trivial, the set of free $\Oc$-modules is
a set of generators.  
The argument that $F_o$ is compact is similar because the map
\[ \bigoplus_{i\in I} \, [F_o, M_i]^{\Ho(\mbox{\scriptsize mod-}\Oc)} \ 
\varrow{1cm} \ [F_o,\coprod_{i\in I} M_i]^{\Ho(\mbox{\scriptsize mod-}\Oc)}  \]
is isomorphic to the map
\[ \bigoplus_{i\in I} \, [S, M_i(o)]^{\Ho(\spec)}\ 
\varrow{1cm} \ [S, \coprod_{i\in I} M_i(o)]^{\Ho(\spec)}\]
and $S$ is compact.   

The proof of (iii) follows as in~\cite[4.3]{ss}. 
The restriction functor $\Psi^*$ preserves objectwise fibrations and
objectwise equivalences, so restriction and extension of
scalars form a Quillen adjoint pair.
For every cofibrant right $\Oc$-module $N$, the induced map  
\[ N \ \iso \ N \sm_{\Oc} \Oc \ \to \ N \sm_{\Oc} \R \]
is an objectwise stable equivalence,
by a similar `cell induction' argument as for ring spectra~\cite[5.4.4]{hss} 
or~\cite[12.7]{mmss}. Thus if $M$ is any right $\R$-module, an
$\Oc$-module map $N \to \Psi^* M$ is an objectwise stable equivalence
if and only if the adjoint $\R$-module map $\Psi_* N = N \sm_{\Oc} \R \to M$ 
is an objectwise stable equivalence.
\end{proof}

\subsection{Quasi-equivalences}
\label{sub-homotopy invariance}
 
This section introduces {\em quasi-equivalences} which
are a bookkeeping device for producing stable equivalences 
between symmetric ring spectra or spectral categories, 
see Lemma \ref{lemma-balanced} below.
The name is taken from \cite[Summary, p.\ 64]{keller-derivingDG},
where this notion is discussed in the context of differential graded algebras.
Every (stable) equivalence of ring spectra gives rise to a quasi-equivalence;
conversely the proof of  Lemma \ref{lemma-balanced} shows
that a single quasi-equivalence encodes a zig-zag of four stable equivalences 
relating two ring spectra or spectral categories.
One place where quasi-equivalences arise `in nature' is the proof 
that weakly equivalent objects in a model category
have weakly equivalent endomorphism monoids, 
see Corollary \ref{cor-ho invariance of endo rings}.

If $\R$ and $\Oc$ are spectral categories, their {\em smash product} 
$\R\sm\Oc$ is the spectral category whose set of 
objects is the cartesian product of the objects of $\R$ and $\Oc$ 
and whose morphism objects are defined by the rule
\[ \R\sm\Oc((r,o),(r',o')) \ = \ \R(r,r') \sm \Oc(o,o') \ . \]
An {\em $\R$-$\Oc$-bimodule} is by definition an 
$\R\rscript{op}\sm\Oc$-module. 
Since modules for us are always contravariant functors, an $\R$-$\Oc$-bimodule
translates to a {\em covariant} spectral functor from 
$\Oc\rscript{op}\sm\R$ to $\spec$.

\begin{definition} \label{def-quasi-equivalence} 
{\em Let $\R$ and $\Oc$ be two spectral categories 
{\em with the same set $I$ of objects}. Then a {\em quasi-equivalence} 
between $\R$ and $\Oc$ is an $\R$-$\Oc$-bimodule 
$M$ together with a collection of `elements' $\varphi_i\in M(i,i)$
(i.e., morphisms $S\to M(i,i)$) for all $i\in I$ 
such that the following holds:
for all pairs $i$ and $j$ of objects the right multiplication with $\varphi_i$
and the left multiplication with $\varphi_j$,
\[ \R(i,j) \ \varr{1cm}{\cdot\varphi_i} \ M(i,j) \ 
\varl{1cm}{\varphi_j\cdot} \ \Oc(i,j) \]
are stable equivalences.}
\end{definition}

\begin{remark} {\em In the important special case of spectral categories 
with a single object, i.e., for two symmetric ring spectra $R$ and $T$, 
a quasi-equivalence is an $R$-$T$-bimodule $M$ together with an element 
$\varphi \in M$ (i.e., a vertex of the 0-th level of $M$ or equivalently
a map $S\to M$ of symmetric spectra) such that the 
left and right multiplication maps with $\varphi$,
\[ R \ \varr{1cm}{\cdot\varphi} \ M \ \varl{1cm}{\varphi\cdot}\ T \]
are stable equivalences of symmetric spectra.}
\end{remark}

If $\Psi:\Oc\to\R$ is a stable equivalence of spectral categories, 
then the target $\R$ becomes an $\R$-$\Oc$-bimodule 
if $\Oc$ acts on the right via $\Psi$.
Furthermore the identity elements in $\R(i,i)$ for all objects $i$ of $\R$
make the bimodule $\R$ into a quasi-equivalence between $\R$ and $\Oc$. 
The following lemma shows conversely that quasi-equivalent 
spectral categories are related by a chain of weak equivalences:

\begin{lemma} \label{lemma-balanced}
Let $\R$ and $\Oc$ be two spectral categories with the same set
of objects. If a quasi-equivalence exists between $\R$ and $\Oc$,
then there is a chain of stable equivalences between $\R$ and~$\Oc$.
\end{lemma}

\begin{proof} 
i) Special case: suppose there exists a quasi-equivalence
$(M,\{\varphi_i\}_{i\in I})$ for which all of the right multiplication maps 
$\cdot\varphi_i :\R(i,j) \to M(i,j)$ are trivial fibrations.
In this case we define a new spectral category $\E(M,\varphi)$ 
with objects $I$ as the pullback of $\R$ and $\Oc$ over $M$.
More precisely for every pair $i,j\in I$ the homomorphism object
$\E(M,\varphi)(i,j)$ is defined as the pullback in $\spec$ of the diagram
\[ \R(i,j) \  \varr{1cm}{\cdot\varphi_i} \ M(i,j) \
\varl{1cm}{\varphi_j\cdot}  \Oc(i,j) \ . \]
Using the universal property of the pullback there is a unique way 
to define composition and identity morphisms in $\E(M,\varphi)$ 
in such a way that the maps $\E(M,\varphi) \to \Oc$ 
and $\E(M,\varphi)\to \R$ are homomorphisms of spectral categories.
  
Since $M$ is a quasi-equivalence, 
all the maps in the defining pullback diagrams are weak equivalences. 
By assumption the horizontal ones are even trivial fibrations, 
so both base change maps $\E(M,\varphi) \to \Oc$ 
and $\E(M,\varphi) \to \R$ are pointwise equivalences 
of spectral categories. The same argument works if instead 
of the right multiplication maps all the left multiplication maps
$\varphi_j\cdot :\Oc(i,j) \to M(i,j)$ are trivial fibrations.

ii) General case: taking fibrant replacement if necessary we can assume 
that the bimodule $M$ is objectwise fibrant. 
The `element'  $\varphi_j$ of $M(j,j)$ corresponds to a map 
$F_j=\Oc(-,j) \to M(-,j)$ from the free $\Oc$-module;
the map is left multiplication by $\varphi_j$ and is 
an objectwise equivalence since $M$ is a quasi-equivalence.
We factor this $\Oc$-module equivalence as a trivial cofibration 
$\alpha_j:F_j \to N_j$ followed by a trivial fibration  
$\psi_j : N_j \to M(-,j)$; in particular, the objects $N_j$ so obtained
are cofibrant and fibrant.
We let $\E(N)$ denote the endomorphism spectral category of the 
cofibrant-fibrant replacements, i.e., the full spectral subcategory 
of the category of $\Oc$-modules with objects $N_j$ for $j\in I$. 
Now we appeal twice to the special case that we already proved, obtaining
a chain of four stable equivalences of spectral categories
\[ \Oc \ \varl{1cm}{\sim} \ \E(W,\alpha) \ \varr{1cm}{\sim} \ \E(N) 
\ \varl{1cm}{\sim} \ \E(V,\psi) \ \varr{1cm}{\sim} \ \R \ . \] 
  
In more detail, we define a $\E(N)$-$\Oc$-bimodule $W$ by the rule
\[ W(i,j) \ = \ \Hom_{\text{mod-}\Oc}(F_i,N_j) \ \iso \ N_j(i) \ . \]
The bimodule $W$ is a quasi-equivalence with respect to the maps $\alpha_j$.
Moreover, the right multiplication map $\cdot\alpha_i$ is the restriction map
\[ \alpha_i^* \ : \ \Hom_{\text{mod-}\Oc} (N_i,N_j) \ \varrow{1cm} \ 
\Hom_{\text{mod-}\Oc}(F_i,N_j) \ . \] 
So $\alpha_i^*$ is a trivial fibration since  $\alpha_i$ 
is a trivial cofibration of $\Oc$-modules and $N_j$ is a fibrant module. 
Case i) above then provides a chain of stable equivalences between $\Oc$ 
and $\E(N)$, passing through $\E(W,\alpha)$.

Now we define an $\R$-$\E(N)$-bimodule $V$ by the rule
\[ V(i,j) \ = \ \Hom_{\text{mod-}\Oc}(N_i,M(-,j)) . \]
The bimodule $V$ is a quasi-equivalence with respect to the maps $\psi_j$.
The left multiplication map $\psi_j\cdot$ is the composition
\[ (\psi_j)_* \ : \ \Hom_{\text{mod-}\Oc} 
(N_i,N_j) \ \varrow{1cm} \ \Hom_{\text{mod-}\Oc}(N_i,M(-,j)) \ . \] 
This time $(\psi_j)_*$ is a trivial fibration since $\psi_j$ is 
a trivial fibration of $\Oc$-modules and $N_i$ is a cofibrant module. 
Furthermore the right multiplication map 
\[ \cdot\psi_i: \R(i,j) \ \varrow{1cm} \
\Hom_{\text{mod-}\Oc}(N_i,M(-,j)) \]
is an equivalence because its composite with the map
\[ \alpha_i^* \ : \ \Hom_{\text{mod-}\Oc}(N_i,M(-,j)) \ \varrow{1cm} \ 
\Hom_{\text{mod-}\Oc}(F_i,M(-,j)) \ \iso \ M(i,j) \]
is right multiplication by $\psi_i$, an equivalence by assumption.
Recall $M$ is objectwise fibrant, so $\alpha_i^*$ is a weak equivalence.
So case i) gives a chain of pointwise equivalences between $\R$ and $\E(N)$, 
passing through $\E(V,\psi)$.
\end{proof}

As a corollary we obtain the homotopy invariance of endomorphism
spectral categories under spectral Quillen equivalences.

\begin{corollary} \label{cor-ho invariance of endo rings}
Suppose $\C$ and $\D$ are  spectral model categories and 
$L:\C \to \D$ is the left adjoint of a  spectral Quillen equivalence.
Suppose $I$ is a set, $\{P_i\}_{i\in I}$ and $\{Q_i\}_{i\in I}$ 
are sets of cofibrant-fibrant objects of $\C$ and $\D$ respectively, and
that for all $i\in I$, $LP_i$ is weakly equivalent to $Q_i$ in $\D$.
Then the spectral endomorphism categories of $\{P_i\}_{i\in I}$ 
and $\{Q_i\}_{i\in I}$ are stably  equivalent.
In particular the spectral endomorphism category of $\{P_i\}_{i\in I}$
depends up to pointwise equivalence only on the weak equivalence
type of the objects $P_i$.
\end{corollary}
\begin{proof} Since the object $LP_i$ is cofibrant and weakly equivalent to
the fibrant object $Q_i$, we can choose a weak equivalence
$\varphi_i:LP_i\to Q_i$ for every $i\in I$. We claim that the collection
of homomorphism objects $\Hom_{\D}(LP_i,Q_j)$ forms a quasi-equivalence 
for the endomorphism spectral categories of  $\{P_i\}_{i\in I}$ 
and $\{Q_i\}_{i\in I}$  with respect to the equivalences $\varphi_i$.
Indeed the endomorphism category of $\{Q_i\}_{i\in I}$ acts on the left by
composition; also right multiplication by $\varphi_j$ is a stable
equivalence since $Q_j$ is fibrant and  $\varphi_j$ is a weak equivalence
between cofibrant objects. 
If $R$ denotes the right adjoint of $L$, then $\Hom_{\D}(LP_i,Q_j)$ is 
isomorphic to $\Hom_{\C}(P_i,RQ_j)$, so the endomorphism category of 
$\{P_i\}_{i\in I}$ acts on the right by composition. Since $R$ and $L$ form
a spectral Quillen equivalence, the adjoints $\widehat{\varphi_i}:P_i\to RQ_i$ 
are weak equivalences between fibrant objects;
so left multiplication by $\varphi_i$ is a stable equivalence 
since $P_i$ is cofibrant. 
The last statement is the special case where $\D=\C$ and $L$ is the identity
functor.
\end{proof}

\section{Eilenberg-Mac Lane spectra and chain complexes}
\label{app-EM}

The proof of the generalized tilting theorem in Section \ref{sec-tilting}
uses the Eilenberg-Mac Lane spectral category $H\uA$ of a ringoid $\uA$.
Recall that a {\em ringoid} is a small category whose hom-sets carry an
abelian group structure for which composition is bilinear. 
A right {\em module} over a ringoid is a contravariant
additive functor to the category of abelian groups.
The Eilenberg-Mac Lane spectral category $H\uA$ of a ringoid $\uA$
is defined in~\ref{def-EM spectral category}.
In this appendix we provide some general facts about 
Eilenberg-Mac Lane spectral categories. 
The main results are that module spectra 
over the Eilenberg-Mac Lane spectral category $H\uA$ are Quillen equivalent to
chain complexes of $\uA$-modules (Theorem \ref{thm-chains and EM}) and
that Eilenberg-Mac Lane spectral categories are determined 
up to stable equivalence by their coefficient ringoid 
(Theorem \ref{prop-EMuniqueness}). 
These properties are not unexpected, and variations have
been proved for the special case of ring spectra in different frameworks.
Indeed the Quillen equivalence of Theorem
\ref{thm-chains and EM} is a generalization and strengthening
of the fact first proved in \cite{robinson-derived} that the unbounded
derived category of modules over a ring $R$ is equivalent to
the homotopy category of $HR$-modules, see also
\cite[IV Thm.\ 2.4]{ekmm} in the context of $S$-algebras.

\subsection{Chain complexes and module spectra}

Throughout this section we fix a ringoid $\uA$, and we want to prove 
Theorem \ref{thm-chains and EM} relating the modules over the 
Eilenberg-Mac Lane spectral category $H\uA$ to complexes of $\uA$-modules
by a chain of Quillen equivalences.
We do not know of a Quillen functor pair which does the job in a single step.
Instead, we compare the two categories through the intermediate
model category of {\em naive $H\uA$-modules}, obtaining a chain
of Quillen equivalences
\vspace{-0cm}
\[
\raisebox{0.18cm}{mod-$H\uA$}
\parbox{2cm}{\begin{center} \raisebox{-0.25cm}{$\scriptstyle  U$} \\ 
\hspace*{0.7cm}
$\begin{diagram} \arrow{e} \end{diagram}$ \\ \vspace{-0.30cm}  
\hspace*{2.2cm} $\begin{diagram} \arrow{w} \end{diagram}$ \\  
\raisebox{0.4cm}{$\scriptstyle L$} \end{center}}
\raisebox{0.18cm}{Nvmod-$H\uA$}
\parbox{2cm}{\begin{center} \raisebox{-0.25cm}{$\scriptstyle \Hc$} \\  
\hspace*{2.2cm}
$\begin{diagram} \arrow{w} \end{diagram}$ \\ \vspace{-0.30cm}  
\hspace*{0.7cm} $\begin{diagram} \arrow{e} \end{diagram}$ \\  
\raisebox{0.4cm}{$\scriptstyle \Lambda$} \end{center}}
\raisebox{0.18cm}{Ch\,$\uA$}
\vspace{-0.5cm}
\]
(the right adjoints are on top), see Corollary~\ref{cor-qeq} and
Theorem~\ref{thm-naive and chains}.  
We mention here that an analogous
statement holds for differential graded modules over a differential
graded ring and modules over the associated Eilenberg-Mac Lane
spectrum, but the proof becomes more complicated; see 
Remark~\ref{rm-alternative EM-model} and~\cite{ss-equiv}. 

\begin{definition} \label{def-naive HA-modules} {\em
Let $\Oc$ be a spectral category.
A {\em naive $\Oc$-module} $M$ consists of a collection $\{M(o)\}_{o\in \Oc}$ 
of $\mN$-graded, pointed simplicial sets together with associative 
and unital action maps 
\[   M(o)_{p} \ \sm \ \Oc(o',o)_q  \ \varrow{1cm} \  M(o')_{p+q} \]
for pairs of objects $o,o'$ in $\Oc$ and for natural numbers $p,q\geq 0$.
A morphism of naive $\Oc$-modules $M\to N$ consists of maps of
graded spaces $M(o) \to N(o)$ strictly compatible with the action of $\Oc$. 
We denote the category of naive $\Oc$-modules by \text{\rm Nvmod-}$\Oc$.}
\end{definition}

Note that a naive module $M$ has {\em no symmetric group action} on $M(o)_n$,
and hence there is no equivariance condition for the action maps.
A naive $\Oc$-module has strictly less structure 
than a genuine $\Oc$-module, so there is a forgetful functor
\[ U \ : \ \text{mod-}\Oc \ \varrow{1cm} \ \text{Nvmod-}\Oc \ . \]
 
The {\em free naive $\Oc$-module} $F_o$ at an object $o\in\Oc$ 
is given by the graded spaces $F_o(o')= \Oc(o',o)$ with action maps
\[  F_o(o')_{p} \ \sm \ \Oc(o'',o')_q \ = \ \Oc(o',o)_{p} \ \sm\ \Oc(o'',o')_q 
\ \varrow{1cm} \ \Oc(o'',o)_{p+q} \ = \ F_o(o'')_{p+q} \]
given by composition in $\Oc$.
In other words, the forgetful functor takes the free, genuine $\Oc$-module 
to the free, naive $\Oc$-module.
The free naive module $F_o$ represents evaluation at the object $o\in\Oc$,
i.e., there is an isomorphism of simplicial sets
\begin{equation} \label{F_o represents}
\map_{\text{\rm Nvmod-}\Oc}(F_o,M) \ \iso \ M(o)_0 \end{equation}
which is natural for naive $\Oc$-modules $M$.

If $M$ is a naive $\Oc$-module, then at every object $o\in\Oc$, $M(o)$ has
an underlying spectrum in the sense of Bousfield-Friedlander \cite[\S 2]{BF}
(except that in \cite{BF}, the suspension coordinates appear on the left,
whereas we get suspension coordinates acting from the right).
Indeed, using the unital structure map $S^1\to  \Oc(o,o)_1$ of the
spectral category $\Oc$, the graded space $M(o)$ gets suspension maps 
via the composite
\[  M(o)_p \ \sm \ S^1 \ \varrow{1cm} \ M(o)_p \ \sm \ \Oc(o,o)_1 \ 
\varrow{1cm} \  M(o)_{p+1} \ . \] 
A morphism of naive $\Oc$-modules $f:M\to N$ is an 
{\em objectwise $\pi_*$-isomorphism} if for all $o\in\Oc$ the map 
$f(o): M(o)\to N(o)$ induces an isomorphism of stable homotopy groups.
The map $f$ is an {\em objectwise stable fibration}
if each $f(o)$ is a stable fibration of spectra in the sense of 
\cite[Thm.~2.3]{BF}).
A morphism  of naive $\Oc$-modules is a {\em cofibration} if it has the
left lifting properties for maps which are objectwise $\pi_*$-isomorphisms 
and objectwise stable fibrations.

\begin{theorem} \label{thm-naive HA-modules}
Let $\uA$ be a ringoid.
\begin{enumerate}
\item[(i)] The category of naive $H\uA$-modules 
with the objectwise $\pi_*$-isomorphisms,
objectwise stable fibrations, and cofibrations is a cofibrantly generated,
simplicial, stable model category. 
\item[(ii)] The collection of free $H\uA$-modules 
$\{F_a\}_{a\in\uA}$ forms a set of compact generators for the
homotopy category of naive $H\uA$-modules.
\item[(iii)] Let $\C$ be a stable model category and consider a Quillen adjoint
functor pair
\vspace{-0cm}
\[
\raisebox{0.18cm}{$\C$}
\parbox{2cm}{\begin{center} \raisebox{-0.25cm}{$\scriptstyle \rho$} \\ 
\hspace*{0.7cm}
$\begin{diagram} \arrow{e} \end{diagram}$ \\ \vspace{-0.30cm}  
\hspace*{2.2cm} $\begin{diagram} \arrow{w} \end{diagram}$ \\  
\raisebox{0.4cm}{$\scriptstyle \lambda$} \end{center}}
\raisebox{0.18cm}{\em Nvmod-$H\uA$}
\]
where $\rho$ is the right adjoint. Then $(\lambda,\rho)$ 
is a Quillen equivalence, provided that
\begin{enumerate}
\item[(a)] for every object $a\in \uA$, the object $\lambda(F_a)$
is fibrant in $\C$ 
\item[(b)] for every object $a\in \uA$, the unit of the adjunction 
$F_a\to \rho\lambda(F_a)$ is an objectwise $\pi_*$-isomorphism, and
\item[(c)] the objects $\{\lambda(F_a)\}_{a\in\uA}$ form a
set of compact generators for the homotopy category of~$\C$.
\end{enumerate}
\end{enumerate}
\end{theorem}
\begin{proof}
(i) We use Theorem \ref{thm-lifting model structure} to establish 
the model category structure.
The category of naive $H\uA$-modules is complete and cocomplete and
every naive $H\uA$-module is small.
The objectwise $\pi_*$-isomorphisms are closed under the 2-out-of-3
condition (Theorem \ref{thm-lifting model structure} (1)).

As generating cofibrations $I$ we use the collection of maps
\[ (\partial\Delta^i)^+ \sm F_a[n]  \ \varr{1cm}{} \
(\Delta^i)^+ \sm  F_a[n]  \]
for all $i,n\geq 0$ and $a\in \uA$.
Here $\Delta^i$ denotes the simplicial $i$-simplex and $\partial\Delta^i$ 
is its boundary; the square bracket $[n]$ means shifting (reindexing) of a
naive $H\uA$-modules and smashing of a module and a pointed simplicial set
is levelwise.
Since the free modules represent evaluation at an object
(see \eqref{F_o represents} above), 
the $I$-injectives are precisely the maps which are objectwise
level acyclic fibrations.

As generating acyclic cofibrations $J$ we use the union
$J=J\rscript{lv}\cup J\rscript{st}$. Here $J\rscript{lv}$ is the set of maps
\[ (\Lambda_k^i)^+ \sm F_a[n]  \ \varr{1cm}{\sim} \
(\Delta^i)^+ \sm F_a[n] \]
for $i,n\geq 0$, $0\leq k\leq i$ and $a\in \uA$,
where $\Lambda^{i,k}$ is the $k$-th horn of the $i$-simplex.
The $J\rscript{lv}$-injectives are the objectwise level fibrations.
Finally, $J\rscript{st}$ consists of the mapping cylinder inclusions 
of the maps
\begin{equation} \label{J^st} \qquad  S^1 \sm \, (\Delta^i)^+ \sm F_a[n{+}1]\
\cup_{ S^1\sm \,(\partial\Delta^i)^+\sm F_a[n{+}1]}\
(\partial\Delta^i)^+ \sm F_a[n]
\ \varrow{1cm} \ (\Delta^i)^+ \sm F_a[n] \ . \end{equation}
Here the mapping cylinders are formed on each simplicial level, 
just as in~\cite[3.1.7]{hss}.
Every map in $J$ is an $I$-cofibration, hence every 
relative $J$-cell complex is too; we claim that in addition, 
every map in $J$ is an objectwise injective $\pi_*$-isomorphism.
Since this property is closed under infinite wedges, pushout, 
sequential colimit and retracts this implies that every relative $J$-cell 
complex is an objectwise injective $\pi_*$-isomorphism and so
condition (2) of Theorem \ref{thm-lifting model structure} holds.

The maps in $J\rscript{lv}$ are even objectwise injective
level-equivalences, so it remains to check the maps in
$J\rscript{st}$. These maps are defined as mapping cylinder inclusions,
so they are injective, and we need only check that the maps in 
\eqref{J^st} above are objectwise $\pi_*$-isomorphisms.
This in turn follows once we know that the maps
\begin{equation} \label{right action map}
 S^1\sm F_a[n{+}1] \ \varrow{1cm} \   F_a[n]  \end{equation}
are objectwise $\pi_*$-isomorphisms.
At level $p\geq n+1$ and an object $b\in\uA$, this map is given 
by the inclusion
\[  S^1\sm (\uA(b,a)\tensor \widetilde{\mathbb Z}[S^{p-n-1}]) \ \varrow{1cm} 
\uA(b,a)\tensor \widetilde{\mathbb Z}[S^{p-n}] \]
whose adjoint is a weak equivalence.
This map is roughly $2(p-n)$-connected, so in the limit we indeed 
obtain a $\pi_*$-isomorphism. 

It remains to check condition (3) of Theorem \ref{thm-lifting model structure},
namely that the $I$-injectives coincide with the maps which are 
both $J$-injective and objectwise $\pi_*$-isomorphisms.
Every map in $J$ is an $I$-cofibration, so $I$-injectives are $J$-injective.
Since $I$-injectives are level acyclic fibrations, they are also
objectwise $\pi_*$-isomorphisms. 
Conversely, suppose $f:M\to N$ is an objectwise $\pi_*$-isomorphism of
naive $H\uA$-modules which is also $J$-injective.
Since $f$ is $J\rscript{lv}$-injective, it is an objectwise level fibration. 
Since $f$ is $J\rscript{lv}$-injective, at every object
$a\in\uA$, the underlying map of spectra $f(a):M(a)\to N(a)$
has the right lifting property for the maps
\[ S^1\sm (\Delta^i)^+  \sm S[n{+}1] 
\cup_{S^1\sm (\partial\Delta^i)^+\sm  S[n{+}1]}\,
(\partial\Delta^i)^+ \sm S[n]
\ \varrow{1cm} \ (\Delta^i)^+ \sm S[n] \ , \]
where $S$ is the sphere spectrum. But then $f(a)$ is a stable fibration
of spectra~\cite[A.3]{sch-theories}, so 
$f$ is an objectwise stable fibration and $\pi_*$-isomorphism. 
By~\cite[A.8 (ii)]{BF}, $f$ is then an objectwise level fibration,
so it is $I$-injective. 
So conditions (1)-(3) of Theorem \ref{thm-lifting model structure}
are satisfied and this theorem provides the model structure.
We omit the verification that the model structure for naive $H\uA$-modules
is simplicial and stable; the latter is a consequence of the
fact that stable equivalences of $H\uA$-modules are defined objectwise and
spectra form a stable model category.

(ii) The stable model structure for naive $H\uA$-modules
is defined so that evaluation at $a\in \uA$ is a right Quillen functor
to the stable model category of Bousfield-Friedlander type spectra.
Moreover, evaluation at $a\in \uA$ has a left adjoint which takes 
the sphere spectrum $S$ to the free module $F_a$.
So the derived adjunction provides an isomorphism of graded abelian groups
\[ [F_a,M]_*^{\text{Ho(Nvmod-}H\uA)}  \ \iso \  
[S,M(a)]_*^{\Ho(\Spc)} \ \iso \ \pi_*\, M(a) \ . \] 
This implies that in the homotopy category
the free modules detect objectwise $\pi_*$-isomorphisms,
so they form a set of generators. It also implies that the  
representable modules are compact, because evaluation at $a\in \uA$
and homotopy groups commute with infinite sums.

(iii) We have to show that the derived adjunction 
on the level of homotopy categories
\vspace{-0cm}
\[
\raisebox{0.18cm}{$\HoC$}
\parbox{2cm}{\begin{center} \raisebox{-0.25cm}{$\scriptstyle R\rho$} \\ 
\hspace*{0.7cm}
$\begin{diagram} \arrow{e} \end{diagram}$ \\ \vspace{-0.30cm}  
\hspace*{2.2cm} $\begin{diagram} \arrow{w} \end{diagram}$ \\  
\raisebox{0.4cm}{$\scriptstyle L\lambda$} \end{center}}
\raisebox{0.18cm}{Ho(Nvmod-$H\uA$)}
\vspace{-0.4cm}
\]
yields equivalences of (homotopy) categories.
The right adjoint $R\rho$ detects isomorphisms:
if $f:X\to Y$ is a morphism in $\HoC$ such that $R\rho(f)$ is
an isomorphism in the homotopy category of naive $H\uA$-modules,
then for every $a\in \uA$, the map $f$ induces an isomorphism on
$[L\lambda(F_a),-]$ by adjointness. Since the objects 
$L\lambda(F_a)$ generate the homotopy category of $\C$, $f$ is
an isomorphism.
It remains to show that the unit of the derived adjunction 
$\eta_M:M\to R\rho (L\lambda(M))$ on the level of homotopy categories 
is an isomorphism for every $H\uA$-module $M$.
For the free $H\uA$-modules $F_a$ 
this follows from assumptions (a) and (b): by (a), $\lambda(F_a)$
is fibrant in $\C$, so the point set level adjunction unit
$F_a\to \rho\lambda(F_a)$ models the derived adjunction unit, then by (b)
$\eta_M$ is an isomorphism.
The composite derived functor $R\rho \circ L\lambda$ is exact; 
the functor $R\rho$ commutes with coproducts (a formal consequence of (ii)),
hence so does  $R\rho \circ L\lambda$ since $L\lambda$ is a left adjoint.
Hence the full subcategory of those $H\uA$-modules $M$ for which 
the derived unit $\eta_M$ is an isomorphism is a localizing subcategory.
Since it also contains the generating  representable modules, it coincides
with the full homotopy category of naive $H\uA$-modules.
\end{proof}

\begin{remark}{\em The reader may wonder why we do not state 
Theorem \ref{thm-naive HA-modules} for a general spectral category $\Oc$.
The reason is that already the analog of part (i), the existence of 
the stable model structure for naive \Omods, can fail 
without some hypothesis on $\Oc$.
The problem can be located: one needs that the analog of the map
\eqref{right action map},
\[ S^1 \sm F_o[n+1] \ \varrow{1cm} \ F_o[n] \] 
which is given by the action of the suspension coordinates from the left,
induces an isomorphism of homotopy groups, taken with respect to suspension 
on the right. But in general, the effects of left and right suspension on
homotopy groups can be related in a complicated way. 
We hope to return to these questions elsewhere.
} 
\end{remark}

As a corollary, we use the criteria in part (iii) of the previous theorem
to establish the Quillen equivalence between the model category of (right) 
$H\uA$-modules of symmetric spectra and the model category 
of (right) naive $H\uA$-modules. 
These criteria are also used to establish the Quillen equivalence
between naive $H\uA$-modules and chain complexes of $\uA$-modules, see
Theorem~\ref{thm-naive and chains} below.

First we recall a general categorical criterion for the existence
of left adjoints. Recall from \cite[Def.~1.1, 1.17]{adamek-rosicky} that
an object $K$ of a category $\C$ is {\em finitely presentable} if
the hom functor $\Hom_{\C}(K,-)$ preserves filtered colimits.
A category $\C$ is called {\em locally finitely presentable} 
if it is cocomplete and there exists a set $A$ 
of finitely presentable objects such that every object  of $\C$ is a filtered
colimit of objects in $A$. The condition `locally finitely presentable'
implies that every object is small in the sense of \cite[2.1.3]{hovey-book}.
For us the point of this definition is that
every functor between locally finitely presentable categories 
which commutes with limits and filtered colimits has a left adjoint
(this is a special case of \cite[1.66]{adamek-rosicky}).
We omit the proof of the following lemma.

\begin{lemma} Let $\uA$ be a ringoid. 
Then the categories of complexes of $\uA$-modules, 
of (genuine) $H\uA$-modules and of naive $H\uA$-modules 
are locally finitely presentable.
\end{lemma}

\begin{corollary}\label{cor-qeq}
The forgetful functor from $H\uA$-modules to naive $H\uA$-modules
is the right adjoint of a Quillen equivalence.
\end{corollary}
\begin{proof}
The forgetful functor $U$ from $H\uA$-modules to naive $H\uA$-modules
preserves limits and filtered colimits. Since source and target category
are locally finitely presentable, $U$ has a left adjoint $L$
by \cite[1.66]{adamek-rosicky}.
The forgetful functor from symmetric spectra to (non-symmetric) spectra
is the right adjoint of a Quillen functor pair, see \cite[4.2.4]{hss}.
So the forgetful functor $U$ from $H\uA$-modules to naive $H\uA$-modules
preserves objectwise stable equivalences and objectwise stable fibrations.
Thus $U$ and $L$ form a Quillen pair, 
and we can apply part (iii) of Theorem \ref{thm-naive HA-modules}.
The left adjoint $L$  sends the naive free modules $F_a$ to
the genuine free modules, so the relevant adjunction unit in condition (b)
is even an {\em isomorphism}.
For every pair of objects $a,b\in \uA$, the symmetric spectrum
$(LF_a)(b)=H\uA(b,a)$ is a symmetric $\Omega$-spectrum, 
hence stably fibrant, which gives condition (a).
The free modules form a set of compact generators
for the homotopy category of genuine $H\uA$-modules, by 
Theorem~\ref{O-modules}, so condition (c) is satisfied.  
\end{proof}

To finish the proof of Theorem \ref{thm-chains and EM}
we now construct a Quillen-equivalence between naive $H\uA$-modules
and complexes of $\uA$-modules.
We define another {\em Eilenberg-Mac Lane functor} 
\[ \Hc \ : \ \mbox{Ch}\uA \ \varrow{1cm} \ \mbox{Nvmod-}H\uA \]
from the category of chain complexes of $\uA$-modules to the category 
of naive modules over $H\uA$. 

For any simplicial set $K$ we denote by $NK$ the normalized  
chain complex of the free simplicial abelian group generated by $K$. 
So $NK$ is a non-negative dimensional chain complex which in dimension $n$ 
is isomorphic to the free abelian group 
on the non-degenerate $n$-simplices of $K$. 
A functor $W$ from the category of chain complexes Ch$_{\mathbb Z}$ to  
the category of simplicial abelian groups is defined by
\[ (WC)_k \ = \ \hom_{\mbox{\scriptsize Ch}_{\mathbb  
Z}}(N\Delta[k],C) \ . \]
For non-negative dimensional complexes, $W$ is just the inverse to the 
normalized chain functor in the Dold-Kan equivalence between simplicial 
abelian groups and non-negative dimensional chain complexes~\cite[1.9]{dold}.
For an arbitrary complex $C$ there is a natural chain map $NWC\to  
C$ which is an isomorphism in positive dimensions and which expresses  
$NWC$ as the $(-1)$-connected cover of $C$.

For a chain complex of abelian groups $C$ we define a graded space 
by the formula
\[ (\Hc C)_n \ = \ W(C[n]) \]
where $C[n]$ denotes the $n$-fold shift suspension of the complex $C$. 
To define the module structure maps we use
the Alexander-Whitney map, see \cite[2.9]{EM-H(pi)} or \cite[29.7]{may-book}. 
This map is a natural, associative and unital  
transformation of simplicial abelian groups
\[ AW \ : \ W(C) \ \tensor \ W(D) \ \varrow{1cm} \ W(C\,\tensor\,  
D) \ . \]
Here the left tensor product is the dimensionwise tensor product of  
simplicial abelian groups, whereas the right one is the tensor  
product of chain complexes.
The Alexander-Whitney map is neither commutative, nor an isomorphism.
By our conventions the $p$-sphere $S^p$ is the $p$-fold smash  
product of the simplicial circle $S^1=\Delta[1]/\partial\Delta[1]$,
so the reduced free abelian group generated by $S^p$ is the $p$-fold tensor  
product of the simplicial abelian group
$\widetilde{\mathbb Z}[S^1]=W({\mathbb Z}[1])$ (where ${\mathbb  
Z}[1]$ is the chain complex which contains a single copy of the  
group $\mathbb Z$ in dimension 1). 
Since the $p$-th space in the Eilenberg-Mac Lane spectrum $H\uA(a,b)$
is given by $H\uA(a,b)_p= \uA(a,b)\tensor \widetilde{\mathbb Z}[S^p]$,
for every chain complex $D$ of $\uA$-modules 
the Alexander-Whitney map gives a map
\begin{eqnarray*}
\Hc(D(b))_p \sm H\uA(a,b)_q & \varrow{1cm} & \Hc(D(b))_p \tensor H\uA(a,b)_q \\
 & \varr{1cm}{\iso} & W(D(b)[p]) \tensor
\uA(a,b)\tensor \underbrace{W({\mathbb Z}[1])\tensor\cdots\tensor W({\mathbb  
Z}[1])}_q) \\
& \varr{1cm}{AW} & W\left(D(b)[p] \tensor \uA(a,b) \tensor 
\underbrace{{\mathbb Z}[1]\tensor\cdots\tensor 
{\mathbb Z}[1]}_q  \right) \\
& \varrow{1cm} &  W(D(a)[p+q]) \ = \  \Hc(D(a))_{p+q}\ .
\end{eqnarray*}
These maps make $\Hc D$ into a naive $H\uA$-module.
The spectra underlying  $\Hc D(a)$ are always $\Omega$-spectra and the  
stable homotopy groups of $\Hc D(a)$ are naturally isomorphic
to the homology groups of the chain complex $D(a)$,
\begin{equation} \label{pi HD = H D} \pi_* \, \Hc D \ \iso \ H_* D 
\end{equation}
as graded $\uA$-modules.

\begin{remark} \label{rm-alternative EM-model}
{\em The functor $\Hc$ should not be confused with
the Eilenberg-Mac Lane functor $H$ 
of Definition \ref{def-EM spectral category}.
The functor $H$ takes values in symmetric spectra, but it cannot
be extended in a reasonable way to chain complexes;
the functor $\Hc$ is defined for complexes, but it only takes
values in {\em naive} $H\uA$-modules. 

The essential difference between the two functors can already be seen
for an abelian group $A$.
The simplicial abelian group $(\Hc A)_n=W(A[n])$ is the minimal model
of an Eilenberg-Mac Lane space of type $K(A,n)$ and it is determined 
by the property that its 
normalized chain complex consists only of one copy of $A$ in dimension $n$. 
The simplicial abelian group $(HA)_n=A\tensor \widetilde{\mathbb Z}[S^n]$ 
is another Eilenberg-Mac Lane space of type $K(A,n)$, but it has
non-degenerate simplices in dimensions smaller than $n$.
The Alexander-Whitney map gives a
weak equivalence of simplicial abelian groups 
$A\tensor \widetilde{Z}[S^n] \to W(A[n])$. 

However, the Alexander-Whitney map is {\em not} commutative, and 
for $n\geq 2$ there is no $\Sigma_n$-action on the minimal model $W(A[n])$ 
which admits an equivariant weak equivalence from 
$\widetilde{Z}[S^n] \tensor A$. More generally, the graded space
$\Hc A$ {\em cannot} be made into a symmetric spectrum which is 
level equivalent to the symmetric spectrum $HA$. 
This explains why the comparison
between $HA$-modules and complexes of $A$-modules has to go
through the category of naive $HA$-modules.}
\end{remark}

\begin{theorem}\label{thm-naive and chains} Let $\uA$ be a ringoid. Then the
Eilenberg-Mac Lane functor $\Hc$ is the right adjoint of
a Quillen equivalence between chain complexes of $\uA$-modules
and naive $H\uA$-modules.
\end{theorem}
\begin{proof} The functor $\Hc$ commutes with limits and filtered colimits, 
and since source and target category of $\Hc$ are locally finitely presentable,
a left adjoint  $\Lambda$ exists by \cite[1.66]{adamek-rosicky}.
The Eilenberg-Mac Lane functor takes values in the category of  
$\Omega$-spectra, which are the (stably) fibrant objects in the category of
naive $H\uA$-modules.
Moreover, it takes objectwise fibrations of chain complexes 
(i.e., epimorphisms) to objectwise level fibrations.
Since level fibrations between $\Omega$-spectra are 
stable fibrations, $\Hc$ preserves fibrations.
Because of the isomorphism labeled \eqref{pi HD = H D},
$\Hc$ takes objectwise quasi-isomorphisms of $\uA$-modules 
to objectwise stable equivalences of $H\uA$-modules, so it also 
preserves acyclic fibrations.
Thus $\Hc$ and $\Lambda$ form a Quillen adjoint functor pair.

Now we apply criterion (iii) of Theorem \ref{thm-naive HA-modules}.
Every chain complex of $\uA$-modules 
is fibrant in the projective model structure, so condition (a) holds.
If we consider the free $\uA$-module $\uA(-,a)$,
as a complex in dimension 0, then the identity 
element in $\uA(a,a)\iso \Hc \uA(-,a)(a)_0$ is represented by a map
of naive $H\uA$-modules $\kappa:F_a\to\Hc (\uA(-,a))$. 
By the adjunction and representability isomorphisms
\[ \hom_{\text{Ch}\uA}(\Lambda(F_a),D) \ \iso \  
\hom_{\text{Nvmod-}H\uA}(F_a,\Hc D) \ \iso \ (\Hc D(a))_0 \ \iso \
\hom_{\text{Ch}\uA}(\uA(-,a),D) \ , \]
so the complexes $\Lambda(F_a)$ and $\uA(-,a)$ represent
the same functor.  Thus, the adjoint of $\kappa$ 
is an isomorphism from $\Lambda(F_a)$ to $\uA(-,a)$. 
The adjunction unit relevant for condition (b) is the map 
$\kappa:F_a\to\Hc (\uA(-,a))\iso \Hc\Lambda(F_a)$.
At an object $b\in\uA$ and in dimension $p$, the map $\kappa$ specializes to
the Alexander-Whitney map
\[ F_a(b)_p \ = \ \uA(b,a)\tensor \widetilde{\mathbb Z}[S^p]
\ \varrow{1cm} \  W(\uA(b,a)[p]) \ = \ \Hc (\uA(b,a))_p \ . \] 
Both sides of this map are Eilenberg-Mac Lane spaces of type
$K(\uA(b,a),p)$, the target being the minimal model.
The map is a weak equivalence, so condition (b) 
of Theorem \ref{thm-naive HA-modules} (iii) holds.
The free modules $\uA(-,a)$ (viewed as a complexes in dimension 0) form a set 
of compact generators for the derived category of $\uA$-modules,
so condition (c) is satisfied. 
\end{proof}

\subsection{Characterization of Eilenberg-Mac Lane spectra}
\label{sub-charac-EM}

In this section we show that Eilenberg-Mac Lane spectral categories 
are determined up to stable equivalence by the property that their homotopy 
groups are concentrated in dimension zero.

\begin{proposition} \label{prop-EMuniqueness}
Let $R$ be a spectral category all of whose morphism spectra are 
stably fibrant and have homotopy groups concentrated in dimension zero. 
Then there exists a natural chain of stable equivalences of 
spectral categories between $R$ and $H\underline{\pi_0R}$.
\end{proposition}

The proposition is a special case of the following statement.
Here we call a stably fibrant spectrum {\em connective} if the negative
dimensional stable homotopy groups vanish.

\begin{lemma} Let $I$ be any set.
There are functors $M$ and $E$ from the category of
spectral categories with object set $I$ to itself and natural transformations
\[ \mbox{\em Id} \ \varr{1cm}{\alpha} \ M \ \varl{1cm}{\beta} \  
E \ \varr{1cm}{\gamma} \ H\underline{\pi_0} \]
with the following properties:
for every spectral category $R$ with connective stably fibrant morphism spectra
the maps $\alpha_R$ and $\beta_R$ are stable equivalences 
and the map $\gamma_R$ induces the canonical isomorphism on 
component ringoids.
\end{lemma}

\begin{proof} The strategy of proof is to transfer the corresponding 
statement from the category of Gamma-rings (where it is easy to prove) 
to the category of symmetric ring spectra and extend it to the 
`multiple object case'. 
We use in a crucial way B{\"o}kstedt's hocolim$_I$ construction
\cite{boekstedt}. 
The functors $M, E$ as well as an intermediate functor $D$  
all arise as lax monoidal  functors 
from the category of symmetric spectra to itself, and the natural maps 
between them  are monoidal transformations.
This implies that when we apply them to the morphism spectra of a spectral
category, then the outcome is again a spectral category in a natural way,
and the transformations assemble into spectral functors.
  
The two functors $M$ and $D$ from the category of symmetric spectra to itself 
are defined in \cite[Sec.\ 3]{shipley-thh}. The $n$-th space of the 
symmetric spectrum $MX$ is defined as the homotopy colimit
\[ (MX)_n \ = \ \hocolim_{{\bf k}\in I} \Omega^ k \mbox{Sing}|X_{k+n}| \ . \]
Here $I$ is a skeleton of the category of finite sets and injections
with objects ${\bf k}=\{0,1,\dots,k\}$; 
for the precise definition and the structure maps making this 
a symmetric spectrum see \cite[Sec.\ 3]{shipley-thh}. 
The map $\alpha : X \to MX$ is induced by the inclusion of $X$ 
into the colimit diagram at $k = 0$. 
In the proof of \cite[Prop.\ 3.1.9]{shipley-thh} it is shown that 
the map $\alpha$ is a stable equivalence (even a level equivalence) for every 
stably fibrant symmetric spectrum $X$.

The $n$-th level of the functor $D$ (for `detection' - 
it detects stable equivalences of symmetric spectra) is defined as
\[ (DX)_n \ =  \ \hocolim_{{\bf k}\in I} \Omega^k \mbox{Sing}|X_k \sm S^n| \ , \]
see \cite[Def.\ 3.1.1]{shipley-thh}.  Also in the proof of 
\cite[Prop.\ 3.1.9]{shipley-thh} a natural map $DX \to MX$ is
constructed which we denote $\beta^1_X$ and 
which is a stable equivalence (even a level equivalence) for every 
stably fibrant spectrum $X$.
  
The symmetric spectrum $DX$ in fact arises from a simplicial functor $QX$. 
The value of $QX$ at a pointed simplicial set $K$ is given by
\[  QX(K) \ = \ \hocolim_{{\bf k}\in I} \Omega^k \mbox{Sing}|X_k \sm K| \ . \]
A simplicial functor $F$ can be evaluated on the simplicial spheres 
to give a symmetric spectrum, which we denote $F(\mS)$. 
In the situation at hand we thus have $DX = QX(\mS)$. If we restrict the
simplicial functor $QX$ to the category $\Gamma\rscript{op}$ 
of finite pointed sets we obtain a $\Gamma$-space \cite{segal,BF} 
denoted $\rho QX$. 
Every $\Gamma$-space can be prolonged to a simplicial functor defined 
on the category of pointed simplicial sets \cite[\S 4]{BF}. 
Prolongation is left adjoint to the restriction functor $\rho$, 
and we denote it by $P$. We then set $EX = (P\rho QX)(\mS)$. 
The unit $P\rho QX \to QX$ of the adjunction between
restriction and prolongation, evaluated at the spheres, 
gives a map of symmetric spectra
\[ \beta^2_X : EX = (P\rho QX)(\mS) \ \varrow{1cm} QX(\mS) = DX \ . \]
We claim that if  $X$ is a connective symmetric $\Omega$-spectrum, then 
$\beta^2_X$ is a level equivalence between connective 
symmetric $\Omega$-spectra.
Indeed, if $X$ is a symmetric $\Omega$-spectrum, it is in particular 
convergent in the sense of~\cite[2.1]{ms}. 
By~\cite[Thm.\ 2.3]{ms} and the remark thereafter, the natural map
$QX(K)\to \Omega(QX(\Sigma K))$ is then a weak equivalence for all pointed
simplicial sets $K$.
This implies (see e.g.~\cite[Lemma 17.9]{mmss}) that $QX$ is a {\em linear}
functor, i.e., that it takes homotopy cocartesian squares to homotopy cartesian
squares. 
In particular, $QX$ converts wedges to products, up to weak equivalence,
and takes values in infinite loop spaces, so the restricted 
$\Gamma$-space $\rho QX$ is very special \cite[p.\ 97]{BF}.
By \cite[Thm.\ 4.2]{BF}, $EX = (P\rho QX)(\mS)$ is 
a connected $\Omega$-spectrum.  Since both $EX$ and $DX$ are 
connected $\Omega$-spectra and the map $\beta^2_X$ is an isomorphism at
level 0, $\beta^2_X$ is in fact a level equivalence. 
The map $\beta_X:EX\to MX$ is defined as the composite of the maps 
$\beta^2_X:EX\to DX$ and $\beta^1_X:DX\to MX$; if $X$ is stably fibrant 
and connective, then both of these are level equivalences, 
hence so is  $\beta_X$.

Every $\Gamma$-space $Y$ has a natural monoidal map $Y \to  H\pi_0Y$ 
to the Eilenberg-Mac Lane $\Gamma$-space 
(\cite[\S 0]{segal}, \cite[Sec.\ 1]{sch-gamma}) 
of its component group which induces the canonical isomorphism on $\pi_0$, 
see \cite[Lemma 1.2]{sch-gamma}. (This map is in fact the unit of 
another monoidal adjunction, namely, between the Eilenberg-Mac Lane 
$\Gamma$-space functor and the $\pi_0$-functor.)
In particular there is such a map of $\Gamma$-spaces 
$\rho QX \to H\pi_0(\rho QX)$. From this we get the map 
\[ \gamma_X : EX = (P\rho QX)(\mS) \ \varrow{1cm} \ 
(PH\pi_0(\rho QX))(\mS) =  H\pi_0 X \]
by prolongation and evaluation of the adjunction unit on spheres.
Whenever X is a stably fibrant, the component groups $\pi_0X$ 
and $\pi_0(\rho QX)$ are isomorphic. The symmetric spectrum associated to
the Eilenberg-Mac Lane $\Gamma$-space by prolongation and then restriction 
to spheres is the Eilenberg-Mac Lane model 
of Definition \ref{def-EM spectral category}.
\end{proof}

\end{appendix}

\end{document}